\documentclass[10pt,reqno]{amsart}
\usepackage{amsthm,amssymb}
\usepackage{mathtools}
\usepackage[a4paper,width=160mm,top=27mm,bottom=27mm]{geometry}
\numberwithin{equation}{section}
\usepackage{todonotes}
\usepackage{color}
\usepackage[colorlinks=true,urlcolor=blue,
citecolor=blue,linkcolor=red,linktocpage,pdfpagelabels,
bookmarksnumbered,bookmarksopen]{hyperref}

\newtheorem{theorem}{Theorem}[section]
\newtheorem{lemma}[theorem]{Lemma}
\newtheorem{definition}[theorem]{Definiton}
\newtheorem{proposition}[theorem]{Proposition}

\theoremstyle{definition}
\newtheorem{remx}[theorem]{Remark}
\newenvironment{remark}
  {\pushQED{\qed}\remx}
  {\popQED\endremx}

\begin{document}

\newcommand{\diver}{\operatorname{div}}
\newcommand{\lin}{\operatorname{Lin}}
\newcommand{\curl}{\operatorname{curl}}
\newcommand{\ran}{\operatorname{Ran}}
\newcommand{\kernel}{\operatorname{Ker}}
\newcommand{\la}{\langle}
\newcommand{\ra}{\rangle}
\newcommand{\N}{\mathbb{N}}
\newcommand{\R}{\mathbb{R}}
\newcommand{\C}{\mathbb{C}}
\newcommand{\T}{\mathbb{T}}

\newcommand{\ld}{\lambda}
\newcommand{\fai}{\varphi}
\newcommand{\0}{0}
\newcommand{\n}{\mathbf{n}}
\newcommand{\uu}{{\boldsymbol{\mathrm{u}}}}
\newcommand{\UU}{{\boldsymbol{\mathrm{U}}}}
\newcommand{\buu}{\bar{{\boldsymbol{\mathrm{u}}}}}
\newcommand{\ten}{\\[4pt]}
\newcommand{\six}{\\[-3pt]}
\newcommand{\vare}{\varepsilon}
\newcommand{\pt}{\partial}
\newcommand{\bg}{\Big}


\newcommand{\mA}{{A}}
\newcommand{\mK}{{Q}}
\newcommand{\mS}{{S}}
\newcommand{\mI}{{I}}

\newcommand{\vu}{ u}
\newcommand{\vz}{ z}
\newcommand{\vv}{ v}
\newcommand{\ve}{ e}
\newcommand{\vw}{ w}
\newcommand{\vf}{ f}
\newcommand{\vh}{ h}
\newcommand{\vp}{ \vec P}
\newcommand{\ang}{{\not\negmedspace\nabla}}
\newcommand{\dxy}{\Delta_{x,y}}
\newcommand{\lxy}{L_{x,y}}
\newcommand{\gnsand}{\mathrm{C}_{\mathrm{GN},3d}}
\newcommand{\wmM}{\widehat{{M}}}
\newcommand{\wmH}{\widehat{{E}}}
\newcommand{\wmI}{\widehat{{I}}}
\newcommand{\wmK}{\widehat{{Q}}}
\newcommand{\wmN}{\widehat{{N}}}
\newcommand{\wm}{\widehat{m}}




\address[Yan-Cheng Lv]{\newline\indent School  of  Mathematics  and  Statistics
	\newline\indent
	Southwest University
	\newline\indent
	Chongqing, 400715 China}
\email{\href{mailto:bestlyc123456@163.com}{bestlyc123456@163.com}}

\address[Gui-Dong~Li]{\newline\indent School  of  Mathematics  and  Statistics
	\newline\indent
	Guizhou University
	\newline\indent
	Guiyang, 550025, China}
\email{\href{mailto:bestdong123@163.com}{bestdong123@163.com}}

\address[Chun-Lei Tang]{\newline\indent School  of  Mathematics  and  Statistics
	\newline\indent
	Southwest University
	\newline\indent
	Chongqing, 400715 China}
\email{\href{mailto:tangcl@swu.edu.cn}{tangcl@swu.edu.cn}}

\thanks{(1) Corresponding author: tangcl@swu.edu.cn }

\thanks{(2) This work was supported by  National Natural Science Foundation of China (No.12561019 and No.12371120)}

\subjclass[2020]{35A15; 35B40; 35J60.}
\title[Normalized ground states of Kirchhoff equation]{Normalized ground states of Kirchhoff equation on waveguide manifolds}
\author{Yan-Cheng Lv, Gui-Dong Li, and Chun-Lei Tang}

\begin{abstract}
We study the   nonlinear Kirchhoff equation  on the waveguide manifolds $\R^d\times\T^n$:
\begin{align*}
	-\left( a+b \int_{\mathbb{R}^d\times\mathbb{T}^n}|\nabla_{x,y}u(x,y)|^2\,dxdy\right) \Delta_{x,y}u+\omega u=|u|^\alpha u,
\end{align*} 
where $a,b>0$, $\omega\in \mathbb{R}$, $\alpha\in(0,2^*)$.
We first establish existence of the normalized ground states when $\alpha\in(0,2^*)$. 
 To precisely characterize the $y$-dependence of ground states, we apply a scaling argument to handle the integral term $ \int_{\mathbb{R}^d \times \mathbb{T}^n} |\nabla_{x,y} u(x,y)|^2  dxdy$. To this end, we construct two distinct modified energy functionals, through which we establish the $y$-independence and $y$-dependence of the ground states in two specific regimes: $(i)$ $ d + n \leq 4$ and $\alpha \in (0, \frac{8}{d+n})$;  $(ii)$ $d = n = 1$ and $\alpha \in (8, \infty)$.
\end{abstract}
\keywords{Kirchhoff equation; waveguide manifolds; normalized ground states.}
\maketitle


\section{Introduction and main results}\label{sec:Introduction and main results}
In this paper, we consider the  nonlinear Kirchhoff equation
\begin{align}\label{eq:1}
	-\left( a+b \int_{\mathbb{R}^d\times\mathbb{T}^n}|\nabla_{x,y}u(x,y)|^2\,dxdy\right) \Delta_{x,y}u+\omega u=|u|^\alpha u, \ \ \ (x,y)\in\mathbb{R}^d \times\mathbb{T}^n,
\end{align}
where $a,b>0$, $d,n\in\mathbb{N}^+$, $\omega\in \mathbb{R}$, $\alpha\in(0,2^*)$ ($2^*=+\infty$ when $d+n<3$, and $2^*=4/(d+n-2)$ when $d+n\geq3$), and $\mathbb{T}=\mathbb{R}/(2\pi\mathbb{Z})$ is
the $2\pi$-periodic torus. Moreover, $\Delta_{ x,y} = \Delta_x +\Delta_y$ 
with  $\Delta_x =\sum_{i=1}^d \partial x_i^2$ is the standard Laplacian on $\mathbb{R}^d$ and  $\Delta_y$ is the Laplace-Beltrami operator on $\mathbb{T}^n$, meanwhile, $\nabla_{ x,y}$  adopts a similar definition. 
The problem \eqref{eq:1} was  proposed by Kirchhoff \cite{kirchhoff1883} as an extension of the classical D'Alembert's wave equations for free vibration of elastic strings. Afterwards, in \cite{Lions1978} Lions introduced an abstract functional analysis framework, such that  problem \eqref{eq:1} and its higher dimension extension have been received much attention  (we can refer to \cite{Alves2001,Alves2009,Arosio1996,Cavalcanti2001,D’Ancona1992}).
Due to the existence of item $\int_{\mathbb{R}^d\times\mathbb{T}^n}|\nabla_{x,y}u(x,y)|^2\,dxdy$, problem \eqref{eq:1} is usually referred to as a non-local problem, which means that problem \eqref{eq:1} is no longer a  pointwise identity, thus leading to mathematical difficulties and providing an interesting research topic.

In order to explore solutions  of
the Kirchhoff equation \eqref{eq:1}, there are two substantially distinct options in terms of the
frequency $\omega$. The first one is to fix the frequency $\omega \in \mathbb{R}$, one can apply the variational method to look for the critical points of the associated energy functional $J_\omega:H_{x,y}^1\to\mathbb{R}$ defined by $J_\omega(u):=E(u)+\frac{\omega}{2}\|u\|_2^2$, where
\begin{align}\label{energy}
	E(u):=\frac{a}{2}\|\nabla_{x,y}u\|_2^2+\frac{b}{4}\|\nabla_{x,y}u\|_2^4-\frac{1}{\alpha+2}\|u\|_{\alpha+2}^{\alpha+2}.
\end{align}
 Alternatively, the study on the normalized solutions of  problem \eqref{eq:1},  while the frequency $\omega$ is unknown and will appear as a Lagrange multiplier. Normalized solutions of problem \eqref{eq:1} can be searched as critical points  of $E(u)$ constrained on $\mathcal{S}(c)$, here
 \begin{align}\label{eq constrain}
 	\mathcal{S}(c):=\left\lbrace u\in H^1(\mathbb{R}^d\times\mathbb{T}^n):\|u\|_{2}^2=\int_{\mathbb{R}^d\times\mathbb{T}^n}|u(x,y)|^2\,dxdy=c\right\rbrace,
 \end{align}
in this situation, the solutions with prescribed $L^2$-norm (i.e. mass $c > 0$ is prescribed). This work is primarily concerned with the existence and properties of normalized solutions to the following problem:
\begin{align}\tag{$\mathcal{P}$}\label{eq:1.1}
	\begin{cases}
		-\left( a+b \int_{\mathbb{R}^d\times\mathbb{T}^n}|\nabla_{x,y}u(x,y)|^2\,dxdy\right) \Delta_{x,y}u+\omega u=|u|^\alpha u, \\
		u\in\mathcal{S}(c).
	\end{cases}
\end{align}

  In recent years, the normalized solutions of problem \eqref{eq:1.1} on Euclidean spaces $\mathbb{R}^d$ have been extensively studied. Notably, Ye has made significant contributions in a series of her works \cite{Ye2015ZAMP,Ye2015MMAS,Ye2016ZAMP}. Specifically, in \cite{Ye2015MMAS}, Ye studied the existence of normalized solutions for $\alpha \in (0, \widehat{2}^*)$ (where $\widehat{2}^* = +\infty$ when $d < 3$ and $\widehat{2}^* = \frac{4}{d-2}$ when $d \geq 3$). In \cite{Ye2015ZAMP}, the author considered the case $\alpha = \frac{8}{d}$ and proved the existence of a mountain-pass-type normalized solution, subsequently, Ye \cite{Ye2016ZAMP} further analyzed the asymptotic behavior of normalized solutions with $\alpha = \frac{8}{d}$. Building on these results, Zeng-Zhang \cite{Zeng2017}   classified  normalized solutions with respect to the exponent $\alpha$ and proved that the minimizer, when it exists, is unique up to translations, while also searching for and proving the uniqueness (up to translations) of mountain pass type normalized solution. Additionally, Qi-Zou \cite{Zou2022} investigated the exact number and explicit expressions of positive normalized solutions. It is worth mentioning that other existence results for normalized solutions of problem \eqref{eq:1.1} on $\mathbb{R}^d$ can be found in \cite{licai2023,He2023,Hu2021} and the references therein.  
 Moreover, some other scholars have studied the existence of normalized solutions of problem \eqref{eq:1.1} in bounded domains $\Omega$ (i.e. $\Omega\subset\mathbb{R}^d$), see \cite{Chang2024JGA,Zhu2021CVPDE}.  
 The research in \cite{Zhu2021CVPDE} encompasses the following special cases:  for  equation \eqref{eq bound} under consideration, a global minimizer exists when 
 $\alpha\in(0, \frac{8}{3})$, whereas only a local minimizer is guaranteed when $\alpha\in ( \frac{8}{3}, 6)$. 
 \begin{align}\label{eq bound}
 	\begin{cases}
 		-\left( 1+b \int_{\Omega}|\nabla u(x)|^2\,dx\right) \Delta u+\omega u=|u|^\alpha u, \ \ \ & \text{in }\Omega\subset\mathbb{R}^3, \\
 		u=0,\ \ \ &\text{in }\partial\Omega,\\
 		\int_{\Omega}|u|^{2}dx=c.
 	\end{cases}
 \end{align}
In addition, there  exist a second solution, namely a mountain-pass type solution, which was obtained by \cite{Chang2024JGA}.
  Moreover, there are many other papers show   the existence of normalized solutions  of problem \eqref{eq bound} ($a=1$ and $b=0$) in the bounded domains, and the interested reader is referred to \cite{Noris2014apde,Pierotti2017CVPDE,Song2023CVPDE} and the references therein.

Besides the domain $\mathbb{R}^d$ or $\Omega\subset\mathbb{R}^d$, the mixed type nature of the domain $\mathbb{R}^d \times \mathcal{M}$ also makes the mathematical analysis on problem \eqref{eq:1.1} rather delicate and challenging.  We noticed the case  $a=1$ and $b=0$ for problem \eqref{eq:1.1}, this is
\begin{align}\label{eq1.4}
	- \Delta_{x,y}u+\omega u=|u|^\alpha u \ \ \ \text{in $\R^d\times \mathcal{M}$,}
\end{align} 
 where $\mathcal{M}$ is a $k$-dimensional compact Riemannian manifold.
 Terracini-Tzvetkov-Visciglia \cite{Tzvetkov2014AnalPDE} studied problem \eqref{eq1.4} in the mass-subcritical regime $\alpha\in(0,\frac{4}{d+k})$, and proved the existence of normalized ground states $u_c$. Another interesting result proved in \cite{Tzvetkov2014AnalPDE}  is the existence of a critical number $c_\star\in(0,\infty)$, such that for any $c\in(0,c_\star)$ the normalized ground states $u_c$ coincide with the ones on $\R^d$, while for $c\in(c_\star,\infty)$ the normalized ground states $u_c$ must have \textbf{non-trivial $y$-dependence}. Subsequently, Pierotti-Verzini-Yu \cite{Yu2025arXiv} obtained local minima by studying problem \eqref{eq1.4}  in the  regime $\alpha\in[\frac{4}{d+k},\min\{\frac{4}{d},\frac{4}{d+k-2}\})$. Furthermore,  when the manifold $\mathcal{M}$ is replaced by $\mathbb{T}$ (here $\mathbb{T}=\mathbb{R}\setminus2\pi\mathbb{Z}$ is
the $2\pi$-periodic torus), by introducing a general framework involving the so-called semivirial-vanishing geometry, Luo \cite{Luo2022arXiv,Luo2023JMPA,Luo2024JGA} studied problem \eqref{eq1.4} in the mass-supercritical regime on $\mathbb{R}^d \times\mathbb{T}$. It is noteworthy that  there are currently many research results the well-posed, long time behavior, existence and stability of the corresponding soliton solutions to problem \eqref{eq1.4}, which can be referred to the literatures \cite{Zhao2021JMAA,Tzvetkov2015FMP,Tzvetkov2016RMI,Zhao2021SIAM}.

Motivated by the above results, the purpose of this paper is therefore
to study the existence and $y$-dependence for normalized ground states of problem \eqref{eq:1.1} on $\mathbb{R}^d \times\mathbb{T}^n$. 
Following  \cite{Cazenave1982CMP}, we consider the   constrained minimization problem on $\mathbb{R}^d \times\mathbb{T}^n $:
\begin{align}\label{minimization}
	m_c:=\inf\left\lbrace E(u):u\in \mathcal{S}(c)\right\rbrace,
\end{align}
where $E(u)$ and $\mathcal{S}(c)$ are defined on \eqref{energy} and \eqref{eq constrain} respectively.
Moreover, in order to formulate our main results, we  also define
\begin{align}
		\widehat{E}(u)&:=\frac{a}{2}\|\nabla_{x}u\|_{L_x^2}^2+\frac{b}{4}\|\nabla_{x}u\|_{L_x^2}^4-\frac{1}{\alpha+2}\|u\|_{L_x^{\alpha+2}}^{\alpha+2},\\
		\widehat{\mathcal{S}}(c)&:=\{ u\in H^1(\mathbb{R}^d):\|u\|_{L_x^2}^2=c\},\\
			\widehat{m}_c&:=\inf\left\lbrace \widehat{E}(u):u\in \widehat{\mathcal{S}}(c)\right\rbrace.\label{hat minimization}
\end{align}
By the classical concentration compactness arguments, our first result concerns the existence of normalized ground state  of $m_c$ for a given mass number $c\in(0,\infty)$ when $\alpha\in(0,8/(d+n))$ with $d+n\leq4$.
\begin{theorem}\label{th:1.1}
Assume that $d+n\leq4$ and $\alpha\in(0,\frac{8}{d+n})$. The following existence results hold:
\begin{itemize}
\item[$(i)$] Suppose that 
 $\alpha\in(0,\min\{\frac{8}{d+n},\frac{4}{d}\})$.
Then  $m_c$ has a minimizer $u_c$ for any $c\in (0,\infty)$.

\item [$(ii)$] Assume  $d>n$ and $\alpha\in[\frac{4}{d},\frac{8}{d+n})$. Define
$$
c_-:=\inf\{c\in(0,\infty):m_c<0\}.
$$
Then $c_-\in(0,\infty)$. Moreover, we have
\begin{itemize}
\item[1)] for any $c\in(0,c_-)$,  $m_c$ has no minimizer;
\item[2)] for $c=c_-$,  $m_{c_-}$ has an minimizer $v$ when $\alpha\in(\frac{4}{d},\frac{8}{d+n})$;
\item[3)] for any $c\in(c_-,\infty)$,  $m_c$ has an minimizer $v$.
\end{itemize}
\end{itemize}
\end{theorem}

Due to the boundedness of $\mathbb{T}$, we derive an interesting fact (see \cite{Luo2022arXiv,Tzvetkov2014AnalPDE}): if $u$ is constant along the $y$-direction, problem \eqref{eq:1.1} reduces to the standard Kirchhoff equation on the Euclidean space $\mathbb{R}^d$. It thus becomes natural to ask whether the minimizer obtained in Theorem \ref{th:1.1}  possess a trivial $y$-dependent component or possess a non-trivial $y$-dependent component. 
Note that, the minimizer is the  ground state solution. We obtain the following result.

\begin{theorem}\label{th:1.2}
Assume that $d+n\leq4$ and $\alpha\in(0,\frac{8}{d+n})$. Then for any $c\in(0,\infty)$, let $u_c$ be a minimizer of $m_c$ deduced from Theorem \ref{th:1.1}, the following statements are true:
\begin{itemize}
\item[$(i)$]
Assume that  $\alpha\in(0,\frac{4}{d+n})$.  Then there exists some $c_{*}\in (0,\infty)$ such that  for any $c\in(0,c_{*})$ it holds $\nabla_y u_c= 0$. 

\item[$(ii)$]
Assume that $\alpha\in(2,\min\{\frac{8}{d+n},\frac{4}{d}\})$. Then there exists some $c_{**}\in(0,\infty)$ such that for any $c\in(c_{**},\infty)$ we have $m_{c}<(2\pi)^n\widehat{m}_{(2\pi)^{-n}c}$ and  $\nabla_y u_c\neq 0$. 
\end{itemize}
\end{theorem}
	\begin{remark}
		To the best of our knowledge, this paper is the first to study the existence of normalized solutions to the Kirchhoff equation and the problem of the trivial or non-trivial $y$-dependent component of the ground state in the context of waveguide manifolds. 
	\end{remark}
	\begin{remark}
		Although the arguments in the proof of Theorem \ref{th:1.2} are very similar to those in the literature \cite{Tzvetkov2014AnalPDE} (see also \cite{Luo2022arXiv}), the additional nonlocal term $\int_{\mathbb{R}^d\times\mathbb{T}^n}|\nabla_{x,y}u(x,y)|^2\,dxdy$ complicates the proof in several places. Moreover, the presence of the non-local term 
		$\int_{\mathbb{R}^d\times\mathbb{T}^n}|\nabla_{x,y}u(x,y)|^2\,dxdy$
		introduces a substantial technical challenge: problem \eqref{eq:1.1} fails to be scale-invariant, which means  the argument in \cite{Tzvetkov2014AnalPDE,Luo2022arXiv} cannot be directly adapted to our setting.
	\end{remark}
Notice that,	in \cite{Tzvetkov2014AnalPDE} (as also presented in \cite{Luo2022arXiv}), the authors considered the modified energy $E^\tau(u)$ ($a=1$ and $b=0$),  which is defined as
	\begin{align*}
		E^{\tau}(u)&:=\frac{a\tau^2}{2}\|\nabla_{y} u\|^2_{2}+\frac{b}{2\tau^{\frac{4}{\alpha}-d}}\|\nabla_{x} u\|^2_{2}\|\nabla_{y} u\|^2_{2}+\frac{b\tau^{2+d-\frac{4}{\alpha}}}{4}\|\nabla_{y} u\|^4_{2}\nonumber\\
		&\quad+\frac{a}{2}\|\nabla_{x} u\|^2_{2}+\frac{b}{4\tau^{\frac{4}{\alpha}+2-d}}\|\nabla_{x} u\|^4_{2}
		-\frac{1}{\alpha+2}\|u\|^{\alpha+2}_{\alpha+2}.
	\end{align*} 
	In \cite{Tzvetkov2014AnalPDE}  the crucial scaling argument is essentially equivalent to demonstrating the following:
	\begin{align}\label{m1l}
		\lim_{\tau\to \infty}\mu_{1,\tau}=(2\pi)^n\widehat{m}_{(2\pi)^{-n}},\quad\lim_{\tau\to 0}\mu_{1,\tau}<(2\pi)^n\widehat{m}_{(2\pi)^{-n}}.
	\end{align}
	Here, the minimization problems $\mu_{c,\tau}$ is defined by 
	\begin{align}
		\mu_{c,\tau}:=\inf\{E^\tau(u):u\in \mathcal{S}(c)\}.\label{mu c tau}
	\end{align}
	Due to the presence of the non-local term $\int_{\mathbb{R}^d\times\mathbb{T}^n}|\nabla_{x,y}u(x,y)|^2\,dxdy$, the problem \eqref{eq:1.1} is no longer scale-invariant, the proof of Theorem \ref{th:1.2} cannot simply rely on showing \eqref{m1l}, and then their method cannot be directly applied. To overcome this obstacle,  we consider the modified energy $E_\lambda(u)$, which is defined as 
	\begin{align*}
			E_{\lambda}(u)&:=\frac{a\lambda^{\frac{4-d}{\alpha-2}+2-\frac{d}{2}}}{2}\|\nabla_{y} u\|^2_{2}+\frac{b\lambda}{2}\|\nabla_{x} u\|^2_{2}\|\nabla_{y} u\|^2_{2}+\frac{b\lambda^2}{4}\|\nabla_{y} u\|^4_{2}\nonumber\\
		&\quad+\frac{a\lambda^{\frac{4-d}{\alpha-2}+1-\frac{d}{2}}}{2}\|\nabla_{x} u\|^2_{2}+\frac{b}{4}\|\nabla_{x} u\|^4_{2}
		-\frac{1}{\alpha+2}\|u\|^{\alpha+2}_{\alpha+2}.
	\end{align*}
	 Analogous to \eqref{m1l}, we establish Theorem \ref{th:1.2} by demonstrating the following:
	\begin{align}\label{tau lambda}
		\lim_{\tau\to \infty}\mu_{1,\tau}=(2\pi)^n\widehat{\mu}_{(2\pi)^{-n},\infty}\quad\lim_{\lambda\to 0}m_{1,\lambda}<(2\pi)^n\widehat{m}_{(2\pi)^{-n}}.
	\end{align}
	Here and later in the proof, the minimization problems are defined by 
	\begin{align}
		m_{c,\ld}&:=\inf\{E_{\ld}(u):u\in \mathcal{S}(c)\},&&	\widehat{m}_{c,\lambda}:=\inf\{\widehat{E}_\ld(u):u\in\widehat{\mathcal{S}}(c)\},\label{m minimizing lambda}\\
			\widehat{\mu}_{c,\tau}&:=\inf\{\widehat{E}^\tau(u):u\in\widehat{\mathcal{S}}(c)\},&&
		\widehat{\mu}_{c,{\infty}}:=\inf\{\widehat{E}^{\infty}(u):u\in\widehat{\mathcal{S}}(c)\},\label{h m c}
	\end{align}
	where
	\begin{align*}
			\widehat{E}_\ld(u)&:=\frac{a\lambda^{\frac{4-d}{\alpha-2}+1-\frac{d}{2}}}{2}\|\nabla_x u\|_{L_x^2}^2+\frac{b}{4}\|\nabla_x u\|_{L_x^2}^4
		-\frac{1}{\alpha+2}\|u\|_{L_x^{\alpha+2}}^{\alpha+2},\\
		\widehat{E}^\tau(u)&:=\frac{a}{2}\|\nabla_x u\|_{L_x^2}^2+\frac{b}{4\tau^{\frac{4}{\alpha}+2-d}}\|\nabla_x u\|_{L_x^2}^4
		-\frac{1}{\alpha+2}\|u\|_{L_x^{\alpha+2}}^{\alpha+2},\\
		\widehat{E}^{\infty}(u)&:=\frac{a}{2}\|\nabla_x u\|_{L_x^2}^2
		-\frac{1}{\alpha+2}\|u\|_{L_x^{\alpha+2}}^{\alpha+2}.
	\end{align*}
Next, we delve into problem \eqref{eq:1.1} for the case where $\alpha\in(\frac{8}{d},2^*)$. To ensure $(\frac{8}{d},2^*)\neq\emptyset$, we will limit our investigation to the situation where $d = n = 1$ (we note that if $d=2$ and $n=1$, at the endpoint case $\alpha=\frac{8}{d}$, the problem \eqref{eq:1.1} is both mass-critical and energy-critical
which requires more techniques, so we  do not consider it here.).  It should be noted that, unlike the mass-subcritical case, in the mass-supercritical case, the energy $E(u)$ is unbounded below on the manifold $\mathcal{S}(c)$, so the method at this time is different from the previous one. As an illustrative example, we consider
\begin{align}\label{eq1.19}
	-\left( a+b\int_{\mathbb{R}^d}|\nabla_x u|^2\,dx\right) \Delta_x u+\omega u=|u|^\alpha u.
\end{align}
%
 Naturaly, one can minimize the functional on a Pohozaev manifold, namely, the following minimization problem
 \begin{align}
 	\widehat{M}_c:=\inf\{\wmH(u):u\in\widehat{V}(c)\},\label{hMc}
 \end{align}
 where 
 \begin{align*}
 	\widehat{V}(c):=\{u\in\widehat{\mathcal{S}}(c):\widehat{P}(u)=0\},
 \end{align*}
 and 
\begin{align*}
	\widehat{P}(u):=a\|\nabla_x u\|_{L_x^2}^2+b\|\nabla_x u\|_{L_x^2}^4-\frac{\alpha d}{2(\alpha+2)}\| u\|_{L_x^{\alpha+2}}^{\alpha+2}=0.
\end{align*}
We can see that $\widehat{E}(u)$ is bounded from below and coercive on  $\widehat{V}(c)$. This idea has been carried out in many cases, see for example \cite{Jean-Luo2013PLMS,Jeanjean1997NA,Ye2015MMAS}. 
 In particular, in the case $a=1$ and $b=0$ of problem \eqref{eq1.19}, we pointed out Jeanjean \cite{Jeanjean1997NA} proved that   $\widehat{M}_c$ is achieved by  $u_c$, which also solves problem \eqref{eq:1.1} with some $\omega_c>0$. We intend to define an analogous minimization problem $M_c$ for problem \eqref{eq:1.1} on $\mathbb{R}\times\T$, following the same line in \cite{Jeanjean1997NA}.   Motivated by \cite{Jeanjean1997NA,Luo2022arXiv,Luo2024JGA}, we consider the following Pohozaev type identity (see Lemma \ref{lemma of poho} below):
\begin{align*}
P(u):=a\|\nabla_x u\|_{2}^2+b\|\nabla_{x,y} u\|_{2}^2\|\nabla_x u\|_{2}^2+b\|\nabla_{x} u\|_{2}^4-\frac{\alpha d}{2(\alpha+2)}\| u\|_{\alpha+2}^{\alpha+2}=0.
\end{align*}
Consequently, one can looking for a  minimizer of  $M_c$,
\begin{align}\label{eq:Mc}
	M_c:=\inf\{E(u):u\in V(c)\},
\end{align}
where
\begin{align*}
	V(c):=\{u\in \mathcal{S}(c):P(u)=0\}.
\end{align*}
We now state the following theorem.
\begin{theorem}\label{th:1.3}
	Assume that  $d = n = 1$ and $\alpha\in(8,\infty)$. For any $c\in(0,\infty)$ the minimization problem $M_c$  has a positive minimizer $u_c$. Moreover, $u_c$ solves problem \eqref{eq:1.1} with some $\omega_c>0$.
\end{theorem}

\begin{remark}
However, when attempting to apply Jeanjean's arguments to the current model, several new issues emerge.

\textbf{Existence of a non-vanishing weak limit.} In the Euclidean case (when considering problem \eqref{eq:1.1} on $\mathbb{R}^d$), a minimizing sequence has a non-vanishing weak limit. In fact, by utilizing the virial-vanishing condition $\widehat{P}(u) = 0$, for a bounded minimizing sequence $\{u_k\}$, we have:
$$ 
0<\delta\leq a\|\nabla_xu_k\|_{L^2(\mathbb{R}^d)}^2\leq\frac{\alpha d}{2(\alpha + 2)}\|u_k\|_{L^{\alpha + 2}(\mathbb{R}^d)}^{\alpha + 2}
\leq C\|\nabla_x u_k\|_{L^2(\mathbb{R}^d)}^{\frac{\alpha d}{2}}\|u_k\|_{L^2(\mathbb{R}^d)}^{\frac{4 - \alpha(d - 2)}{2}}
\leq C\|\nabla_x u_k\|_{L^2(\mathbb{R}^d)}^{\frac{\alpha d}{2}},
$$
This, in turn, implies:
$$ \liminf_{k\to\infty}\|u_k\|_{L^{\alpha + 2}(\mathbb{R}^d)}^{\alpha + 2}\geq C
\liminf_{k\to\infty}\|\nabla_x u_k\|_{L^{2}(\mathbb{R}^d)}^{2}>0.$$
Consequently, a non-vanishing weak limit exists. Generally speaking, on $\mathbb{R}\times \T$, merely relying on the  Gagliardo-Nirenberg inequality  is insufficient to obtain a non-vanishing weak limit of a minimizing sequence. We will overcome this difficulty by making use of a generalization of the Gagliardo-Nirenberg inequality \eqref{eq:gn2}, which was presented in \cite{Luo2022arXiv}.

 \textbf{The sign of  $\omega$ in problem \eqref{eq:1.1}.}  The positive frequency parameter $\omega > 0$ plays a crucial role in proving the $L^2$ convergence of the Palais-Smale sequence. In fact, from the Pohozaev identity of problem \eqref{eq:1.1} (refer to Lemma \ref{lemma of poho}), we can derive:
\begin{align*}
	a\|\partial_y u\|_{2}^2 + b\|\nabla_x u\|_{2}^2\|\partial_y u\|_{2}^2 + b\|\partial_y u\|_{2}^4 + \omega\|u\|_{2}^2=\frac{2\alpha + (4 - \alpha d)}{2(\alpha + 2)}\|u\|_{\alpha + 2}^{\alpha + 2}.
\end{align*}
In the Euclidean case,  $\partial_y u = 0$. Thus, the positivity of $\omega$ can be deduced from $u\neq 0$, which is valid for the non-vanishing weak limit we obtained previously. In our case, we cannot rule out the possibility by the Pohozaev identity only, and furthermore we need to analyze the limit equation to obtain the positivity of $\omega$.
\end{remark}

Inspired by Theorem \ref{th:1.2}, we also consider that the $y$-dependence of $u_c$ as given in Theorem \ref*{th:1.3}.

\begin{theorem}\label{th:1.4}
Assume that  $d = n = 1$ and $\alpha\in(8,\infty)$, for any $c\in(0,\infty)$ let $u_c$ be a minimizer of $M_c$ deduced from Theorem \ref{th:1.3}.	 Then there exists some $0<c^*\leq c^{**}$ such that
	\begin{itemize}
		\item [$(i)$] for all $c \in (0,c^*)$ we have $M_{c}<2\pi \widehat{M}_{(2\pi)^{-1}c}$, and any minimizer $u_c$ of $M_{c}$ satisfies $\pt_y u_c\neq 0$;
		
		\item [$(ii)$]   for $c\in(c^{**},\infty)$ any minimizer $u_c$ of $M_{c}$ satisfies $\pt_y u_c=0$.
	\end{itemize}
\end{theorem}
\begin{remark}
	We observe that, owing to the  term $\int_{\mathbb{R}^d\times\mathbb{T}^n}|\nabla_{x,y}u(x,y)|^2\,dxdy$, we are unable to make a comparison between the rescaled variational problem and the unparameterized constant variational problem, similar to what was done in \cite{Luo2022arXiv,Tzvetkov2014AnalPDE}. Consequently, the question of whether $c_* = c_{**}$ (or $c^* = c^{**}$) remains an intriguing and challenging open problem. Here, $c_*$ and $c_{**}$ ($c^*$ and $c^{**}$) are defined in Theorem \ref{th:1.2} (Theorem \ref{th:1.4}). 
\end{remark}
 Similar as the proof of Theorem \ref{th:1.2}, 
proving Theorem \ref{th:1.4}
is essentially equivalent to showing that
\begin{align*}
	\lim_{\ld\to 0}M_{1,\ld}<(2\pi)^n\widehat{M}_{(2\pi)^{-n},\lambda},\quad\lim_{\tau\to \infty}\nu_{1,\tau}=(2\pi)^n\hat{\nu}_{(2\pi)^{-n},\infty}.
\end{align*}
And we need  consider the following minimization problems:
\begin{align}
	M_{1,\ld}&:=\inf\{E_\ld(u):u\in V_\lambda(1)\},&&\widehat{M}_{c,\lambda}:=\inf\{\widehat{E}_\ld(u):u\in\widehat{V}_\lambda(c)\},\label{M c lambda}\\
	\nu_{1,\tau}&:=\inf\{E^\tau(u):u\in V^\tau(1)\},&&\widehat{\nu}_{c,\tau}:=\inf\{\widehat{E}^{\tau}(u):u\in\widehat{V}^\tau(c)\}\label{nu c tau},\\
	\widehat{\nu}_{c,\infty}&:=\inf\{\widehat{E}^{\infty}(u):u\in\widehat{V}^\infty(c)\},\label{hat nu infty},
\end{align}
where
\begin{align*}
	V_\lambda(c)&:=\{u\in \mathcal{S}(c):P_\lambda(u)=0\},\quad&&\widehat{V}_\lambda(c):=\{u\in \widehat{\mathcal{S}}(c):\widehat{P}_\lambda(u)=0\},\\
	V^\tau(c)&:=\{u\in \mathcal{S}(c):P^\tau(u)=0\},\quad&&\widehat{V}^\tau(c):=\{u\in \widehat{\mathcal{S}}(c):\widehat{P}^\tau(u)=0\},\\
	\widehat{V}^\infty(c)&:=\{u\in \widehat{\mathcal{S}}(c):\widehat{P}^\infty(u)=0\}\nonumber,
\end{align*}
and
\begin{align}	
	P_\lambda(u)&:=a\lambda^{\frac{4-d}{\alpha-2}+1-\frac{d}{2}}\|\nabla_x u\|_{2}^2+b\lambda\|\nabla_{x,y} u\|_{2}^2\|\nabla_x u\|_{2}^2+b\|\nabla_x u\|_{2}^4-\frac{\alpha d}{2(\alpha+2)}\| u\|_{\alpha+2}^{\alpha+2}=0,\nonumber\\ 
	P^\tau(u)&:=a\|\nabla_x u\|_{2}^2+\dfrac{b}{\tau^{\frac{4}{\alpha}-d}}\|\nabla_{x,y} u\|_{2}^2\|\nabla_x u\|_{2}^2+\dfrac{b}{\tau^{\frac{4}{\alpha}+2-d}}\|\nabla_x u\|_{2}^4-\frac{\alpha d}{2(\alpha+2)}\| u\|_{\alpha+2}^{\alpha+2}=0, \label{pohozave tau},\\
	\widehat{P}_\lambda(u)&:=a\lambda^{\frac{4-d}{\alpha-2}+1-\frac{d}{2}}\|\nabla_x u\|_{L_x^2}^2+b\|\nabla_x u\|_{L_x^2}^4-\frac{\alpha d}{2(\alpha+2)}\| u\|_{L_x^{\alpha+2}}^{\alpha+2}=0\nonumber\\,
	\widehat{P}^\tau(u)&:=a\|\nabla_x u\|_{L_x^2}^2+\dfrac{b}{\tau^{\frac{4}{\alpha}+2-d}}\|\nabla_x u\|_{L_x^2}^4-\frac{\alpha d}{2(\alpha+2)}\| u\|_{L_x^{\alpha+2}}^{\alpha+2}=0,\nonumber\\
	\widehat{P}^\infty(u)&:=a\|\nabla_x u\|_{L_x^2}^2-\frac{\alpha d}{2(\alpha+2)}\| u\|_{L_x^{\alpha+2}}^{\alpha+2}=0.\label{hat pohozave infty}
\end{align}

For simplicity, we ignore in most cases the dependence of the function spaces on their underlying domains and hide this dependence in their indices. For example $L_x^p=L^p(\mathbb{R}^d)$, $L_y^p=L^p(\mathbb{T}^n)$, $L_{x,y}^p=L^p(\mathbb{R}^d\times \mathbb{T}^n)$ for $p\geq1$, $H_{x,y}^1= H^1(\mathbb{R}^d\times \mathbb{T}^n)$ and so on. The norm $\|\cdot\|_p$ is defined by $\|\cdot\|_{L^p_{x,y}}$.  We denote by "$\to$"
(resp., by "$\rightharpoonup$") the strong convergence (resp., the weak convergence) in corresponding space, and denote by $\hookrightarrow \hookrightarrow $ the compact embedding. $C_1,C_2,\cdots$ denote positive constants possibly different in different places. For our purpose we will also make use of the following sets:
\begin{align}
	C_{\rm per}^\infty(\mathbb{T}^n)&:=\left\lbrace u\in C^\infty(\mathbb{R}^n),\,\forall\,y\in\mathbb{R}^n:u(y+(2\pi)^n)=u(y)\right\rbrace ,\nonumber\\
	C_0^\infty(\mathbb{R}^d)\otimes C_{\rm per}^\infty(\mathbb{T}^n)&:=\left\lbrace \sum_{i=1}^k f_ig_i:k\in\N,f_i\in C_0^\infty(\mathbb{R}^d),g_i\in C_{\rm per}^\infty(\mathbb{T}^n) \right\rbrace .\label{eq1.14}
\end{align}

This paper is organized as follows: in Section 2, we introduce some preliminary knowledge, which will be useful later. In Section 3, we will consider mass subcritical case and  prove Theorems \ref{th:1.1} and \ref{th:1.2} in Subsection \ref{sub1.1} and in Subsection \ref{sub1.2} respectively. In Section 4, we will consider mass supercritical case and prove Theorems \ref{th:1.3} and \ref{th:1.4} in Subsection \ref{sub2.1} and in Subsection \ref{sub2.2} respectively.

\section{Some preliminaries}
Before giving the proofs of our main results, we first collect some useful auxiliary lemmas. The first two types of inhomogeneous Gagliardo-Nirenberg inequality on $\mathbb{R}^d\times\mathbb{T}^n$ (see \cite{Luo2022arXiv,Tzvetkov2014AnalPDE}). 

\begin{lemma}\label{lem:gn1}
On $\mathbb{R}^d\times\mathbb{T}^n$, if $\alpha\in(0,2^*)$ there exists some $C(\alpha,d,n)>0$ such that for all $u\in H_{x,y}^1$ we have
\begin{align}\label{eq:gn1}
\|u\|_{\alpha+2}^{\alpha+2}\leq C(\alpha,d,n) \|u\|^{\frac{\alpha(d+n)}{2}}_{H^1_{x,y}}\|u\|_2^{\alpha+2-\frac{\alpha(d+n)}{2}}.
\end{align}
\end{lemma}
\begin{lemma}\label{lem:gn2}
	On $\mathbb{R}^d\times\T$, if $\alpha\in(0,2^*)$ there exists some $C(\alpha,d)>0$ such that for all $u\in H_{x,y}^1$ we have
	\begin{align}\label{eq:gn2}
		\|u\|_{\alpha+2}^{\alpha+2}\leq C(\alpha,d)\|\nabla_x u\|_2^{\frac{\alpha d}{2}}\|u\|_2^{\frac{4-\alpha(d-1)}{2}}
		\left( \| u\|_{2}^{\frac{\alpha}{2}}+\|\pt_y u\|_{2}^{\frac{\alpha}{2}}\right) .
	\end{align}
\end{lemma}
We also record a very useful lemma from \cite{Luo2022arXiv} (see also \cite[The proof of Theorem 1.1]{Tzvetkov2014AnalPDE}), which plays an important role for identifying a non-vanishing weak limit of a minimizing sequence. 
\begin{lemma}\label{lem:non-vanishing}
	Let $\{u_k\}$ be a bounded sequence in $H_{x,y}^1$. Assume also that 
	\begin{align*}
		\liminf_{n\to\infty}\|u_k\|_{\alpha+2}>0 \quad\text{for $\alpha\in(0,2^*)$}.
	\end{align*}
	Then there exist $\{x_k\}\subset \R^d$ and some $u\in H_{x,y}^1\setminus\{0\}$ such that up to a subsequence
	\begin{align*}
		u_k(x+x_k,y)\rightharpoonup u(x,y)\quad\text{ in $H_{x,y}^1$}.
	\end{align*}
\end{lemma}
\begin{lemma}\label{lemma of poho}
	Let $u\in H_{x,y}^1$ be a solution of problem \eqref{eq:1.1} with some $\omega>0$. Then $P(u)=0$.
\end{lemma}

\begin{proof}
	Using standard elliptic regularity theory (see  \cite[Lemma B.3]{Book1996Struwe}) and bootstrap arguments (see  \cite[Theorem 8.1.1]{BookCazenave2003}) we know that $u\in W_{\rm loc}^{3,q}(\R^{d+n})$ for all $q\in[1,\infty)$. Follow  \cite[Theorem 8.1.2]{BookCazenave2003},  direct computation yields
	\begin{align}
	-K\Delta_x u(x\cdot\nabla_x u)&=-K\frac{d-2}{2}|\nabla_x u|^2+K\nabla_x\cdot
		\left(  \nabla_x u(x\cdot\nabla_x u)+x\frac{|\nabla_x u|^2}{2}\right) ,\label{pohozaev1}\\
	(\omega u-|u|^\alpha u)x\cdot\nabla_xu&=-d\left( \frac{\omega}{2}|u|^2-\frac{1}{\alpha+2}|u|^{\alpha+2}\right) +div_x\left( x\left( \frac{\omega}{2}|u|^2-\frac{1}{\alpha+2}|u|^{\alpha+2}\right) \right) 
		,\label{pohozaev2}\\
		-K\Delta_y u (x\cdot\nabla_x u)&=-K\frac{d|\nabla_y u|^2}{2}+K\nabla_x\cdot\left( x\frac{|\nabla_y u|^2}{2}\right)  -Kdiv_y \left( \nabla_yu (x\cdot\nabla_x u)\right) 
,\label{pohozaev3}
	\end{align}
	where $K:= a+b\int_{\mathbb{R}^d\times\mathbb{T}^n}|\nabla_{x,y}u|^2\,dxdy$. Define
	\begin{align*}
		R_k:=B_k^x\times\T^n,\quad B_k^x=\{x\in\R^d:|x|<k\}.
	\end{align*}
       Integrate \eqref{pohozaev1}-\eqref{pohozaev3} over $R_k$ and followed by divergence theorem and the periodic boundary conditions, we obtain
	\begin{align}
		0=&-\left( K\int_{\T^n}\int_{\pt B_k^x}(\nabla_x u \cdot n_x)(\sigma_x\cdot\nabla_x u)-(\sigma_x\cdot \theta_x)\frac{|\nabla_x u|^2}{2}
		\,d\sigma_xdy\right) \nonumber\\
		&+\int_{\T^n}\int_{\pt B_k^x}(\sigma_x\cdot \theta_x)\left(\frac{\omega}{2}|u|^2-\frac{1}{\alpha+2}|u|^{\alpha+2}\right)
		\,d\sigma_xdy\nonumber\\
		&+\int_{\T^n}\int_{\pt B_k^x}K(\sigma_x\cdot \theta_x)\frac{|\nabla_y u|^2}{2}
		\,d\sigma_xdy\nonumber\\
		&-K\frac{d-2}{2}\|\nabla_x u\|_{L^2(R_k)}^2-K\frac{d}{2}\|\nabla_y u\|_{L^2(R_k)}^2-\frac{d\omega}{2}\|u\|_{L^2(R_k)}^2+\frac{d}{\alpha+2}\|u\|_{L^{\alpha+2}(R_k)}^{\alpha+2}\nonumber\\
		=&I+II+III-K\frac{d-2}{2}\|\nabla_x u\|_{L^2(R_k)}^2-K\frac{d}{2}\|\nabla_y u\|_{L^2(R_k)}^2-\frac{d\omega}{2}\|u\|_{L^2(R_k)}^2+\frac{d}{\alpha+2}\|u\|_{L^{\alpha+2}(R_k)}^{\alpha+2}\nonumber.
	\end{align}
	where $\theta_x = \frac{x}{k}$ denotes the unit outward normal at $x \in\pt B_k^x$. Particularly, the terms $I, II, III$ are well-defined, thanks to the smoothness of $u$. Arguing as in the proof of \cite[Proporsition 1]{Lions1983iARMA}, there exists a sequence $\{k\}\subset(0,\infty)$ with $n\to\infty$ such that
	$$ 
	I+II+III\to 0\quad\text{as}\quad k\to\infty.
	$$
	 Combining with the dominated convergence theorem, we infer that
	\begin{align}\label{1poho}
		\mathcal{P}(u):=K\frac{d-2}{2}\|\nabla_x u\|_{2}^2+K\frac{d}{2}\|\nabla_y u\|_{2}^2+\frac{d\omega}{2}\|u\|_{2}^2-\frac{d}{\alpha+2}\|u\|_{\alpha+2}^{\alpha+2}=0.
	\end{align}
	On the other hand, multiplying problem \eqref{eq:1.1} with $\bar{u}$ and integrating over $\R^d\times\T^n $, we obtain
	\begin{align}\label{2poho}
		0=a\|\nabla_{x,y} u\|_{2}^2+b\|\nabla_{x,y} u\|_{2}^4+\omega\|u\|_{2}^2-\|u\|_{\alpha+2}^{\alpha+2}.
	\end{align}
	Eliminating $\omega\|u\|_2^2$ in \eqref{1poho} and \eqref{2poho} we conclude that $P(u)=0$.
\end{proof}

\section{Existence and $y$-dependence of normalized solution: Mass subcritical case}
\subsection{Existence  of normalized solutions}\label{sub1.1}
In this subsection, we study the existence of normalized solutions of problem \eqref{eq:1.1} in the mass subcritical. 
First, we give some properties of the minimization problems $\widehat{m}_c$ and $\widehat{M}_c$ which are defined in \eqref{hat minimization} and \eqref{hMc} respectively.
\begin{proposition}[Properties of $\widehat{m}_c$ and $\widehat{M}_c$, \cite{Ye2015MMAS}]\label{thm:2.4}
	Assume $d\leq3$. The following statements hold true:
	\begin{itemize}
		\item[(i)] Assume that 	$\alpha\in(0,\frac{8}{d})$. Then $\widehat{m}_c\leq0$ for each $c>0$.  Moreover, there exists some $c_\alpha \geq0$ such that $\widehat{m}_c< 0$ for $c > c_\alpha$ if $\alpha\in(0,\frac{8}{d})$. In particular, $c_\alpha=0$ if $\alpha\in(0,\frac{4}{d})$.
		\item[(ii)] Assume that  $\alpha\in(0,\frac{4}{d}]$.
		Then for any $c\in(c_\alpha,\infty)$ we have $\widehat{m}_c\in(-\infty,0)$ and $\widehat{m}_c$ has a minimizer $u_c\in \widehat{\mathcal{S}}(c)$.
		
		\item[(iii)] Assume  $\alpha\in(\frac{4}{d},\frac{8}{d})$. Then 
		\begin{itemize}
			\item[(a)]For any $c\in(0,c_{\alpha})$ we have $\widehat{m}_c=0$ and $\widehat{m}_c$ has no minimizer.
			\item[(b)]For any $c\in[c_{\alpha},\infty)$ we have $\widehat{m}_c\in(-\infty,0)$ and $\widehat{m}_c$ has a minimizer $u_c\in \widehat{\mathcal{S}}(c)$.
		\end{itemize}
		\item[(iv)] Assume that  $\alpha\in(\frac{8}{d}, \widehat{2}^*)$.
		Then for any $c\in(0,\infty)$ we have $\widehat{M}_c$ has a minimizer $u_c\in \widehat{\mathcal{S}}(c)$.
	\end{itemize}
\end{proposition}

 \begin{proposition}\label{cor 3.1}
 	When $\alpha\in(0,\frac{8}{d})$, we have $m_c \leq 0$, where $m_c$ is defined in \eqref{minimization}. In particular, $m_c  < 0$ for all $c>0$ when $\alpha\in(0,\frac{4}{d})$.
 \end{proposition}
 \begin{proof}
 		For the following problem
 	\begin{align}\label{eq tilde}
 		\begin{cases}
 			-\left( a+\widetilde{b} \int_{\mathbb{R}^d}|\nabla u(x)|^2\,dx\right) \Delta u+\omega u=|u|^\alpha u, \ \ \ & \text{in }\mathbb{R}^d, \\
 			\int_{\mathbb{R}^d}|u|^{2}dx=c,
 		\end{cases}
 	\end{align}
 	we consider the following constrained minimization problem:
 	\begin{align}\label{tilde minimization}
 		\widetilde{m}_c:=\inf\left\lbrace \widetilde{E}(u):u\in \widehat{\mathcal{S}}(c)\right\rbrace,
 	\end{align}
 	where $\widetilde{E}(u)$ is the energy function of problem \eqref{eq tilde}. Then $\widetilde{m}_c$ also has the properties of $\widehat{m}_c$ as stated in Proposition \ref{thm:2.4}.
 	Now	we assuming that $\widetilde{b}=b(2\pi)^n$ and $u$ is a minimizer of $\widetilde{m}_{(2\pi)^{-n}c}$ in $\widehat{\mathcal{S}}((2\pi)^{-n}c)$ (i.e. $\widetilde{E}(u)=\widetilde{m}_{(2\pi)^{-n}c}$), this is $u\in\mathcal{S}(c)$, then we already have 
 	\begin{align*}
 		m_c\leq E(u)&=\frac{a}{2}\|\nabla_{x,y}u\|_2^2+\frac{b}{4}\|\nabla_{x,y}u\|_2^4-\frac{1}{\alpha+2}\|u\|_{\alpha+2}^{\alpha+2}\\
 		&=(2\pi)^n\left( \frac{a}{2}\|\nabla_{x,y}u\|_{L_x^2}^2+\frac{b(2\pi)^n}{4}\|\nabla_{x,y}u\|_{L_x^2}^4-\frac{1}{\alpha+2}\|u\|_{L_x^{\alpha+2}}^{\alpha+2}\right) \\
 		&=(2\pi)^n\widetilde{E}(u)=(2\pi)^{n}\widetilde{m}_{(2\pi)^{-n}c}.
 	\end{align*}
 	Moreover, since $\widetilde{m}_c<0$ for all $c>0$ when $\alpha\in(0,\frac{4}{d})$, we have  $m_c \leq (2\pi )^{n} \widetilde{m}_{(2\pi )^{-n}c} < 0$ for all $c>0$ when $\alpha\in(0,\frac{4}{d})$. 
 \end{proof}

\begin{lemma}\label{lem:3.1}
	Assume that 	$\alpha\in(0,\frac{8}{d+n})$. Then any minimizing sequence of $m_c$ is bounded in $H^1_{x,y}$. In addition, the mapping $c\mapsto m_c$ is continuous on $(0,\infty)$ and the mapping $c\mapsto c^{-1}m_c$ is  decreasing on $(0,\infty)$. 
	\begin{proof}
		Let $\{u_k\}\subset \mathcal{S}(c)$ be a minimizing sequence of $m_c$.  By using Lemma \ref{lem:gn1} we infer that
		\begin{align}
			&m_c+\frac{ac}{2}-\frac{bc^2}{4}+o_k(1)\nonumber\\
			&= \frac{ac}{2}-\frac{bc^2}{4}+\frac{a}{2}\|\nabla_{x,y}u_k\|_2^2+\frac{b}{4}\|\nabla_{x,y}u_k\|_2^4-\frac{1}{\alpha+2}\|u_k\|_{\alpha+2}^{\alpha+2}\nonumber\\
			&= \frac{ac}{2}-\frac{bc^2}{2}+\left( \frac{a}{2}-\frac{bc}{2}\right) \|\nabla_{x,y} u_k\|_2^2+\frac{b}{4}\|\nabla_{x,y}u_k\|_2^4+\frac{bc}{2} \|\nabla_{x,y} u_k\|_2^2+\frac{bc^2}{4}-\frac{1}{\alpha+2}\|u_k\|_{\alpha+2}^{\alpha+2}\nonumber\\
			&\geq \left( \frac{a}{2}-\frac{bc}{2}\right)\|u_k\|_{H_{x,y}^1}^2+\frac{b}{4}\|u_k\|_{H_{x,y}^1}^4-Cc^{\frac{\alpha+2}{2}-\frac{\alpha(d+n)}{2}}\|u_k\|^{\frac{\alpha(d+n)}{2}}_{H_{x,y}^1}\nonumber\\
			&\geq \inf_{t>0}\left[ \left( \frac{a}{2}-\frac{bc}{2}\right)t^2+ \frac{b}{4}t^4-Cc^{\frac{\alpha+2}{2}-\frac{\alpha(d+n)}{4}}t^{\frac{\alpha(d+n)}{2}}\right],\label{eq:3.1}
		\end{align}
	     this in turn implies that $m_c>-\infty$ for $\alpha<8/(d+n)$. 
		If $\|u_k\|_{H^1_{x,y}}\to \infty$,  using  \eqref{eq:3.1} we deduce the contradiction $m_c\geq \infty$,  so $\{u_k\}$ is bounded in $H_{x,y}^1$. 
		
		New, we consider the continuous of $c\mapsto m_c$.
		Let $\{u_k\}\subset \mathcal{S}(c)$ be a  minimizing sequence of $m_c$ and $\{c_j\}$ be a positive sequence with $\lim_{j\to\infty}c_j=c$. Then by definition of $m_c$ we infer that
		\begin{align}
			m_{c_j}\leq E\left( \frac{\sqrt{c_j}}{\sqrt{c}}u_k\right) &=\frac{ac_j}{2c}\|\nabla_{x,y} u_k\|_2^2+\frac{bc_j^2}{4c^2}\|\nabla_{x,y}u_k\|_2^4
			-\frac{1}{\alpha+2}\bg(\frac{c_j}{c}\bg)^{\frac{\alpha}{2}+1}\|u_k\|_{\alpha+2}^{\alpha+2} \nonumber\\
			&=\frac{c_j}{c}\left( \frac{a}{2}\|\nabla_{x,y} u_k\|_2^2+\frac{b}{4}\|\nabla_{x,y}u_k\|_2^4
			-\frac{1}{\alpha+2}\|u_k\|_{\alpha+2}^{\alpha+2}\right)+\frac{bc_j}{4c}\left( \frac{c_j}{c}-1\right) \|\nabla_{x,y}u_k\|_2^4 
			\nonumber\\
			&\quad+\frac{c_j}{c}\left( 1-\bg(\frac{c_j}{c}\bg)^{\frac{\alpha}{2}}\right) \frac{1}{\alpha+2}\|u_k\|_{\alpha+2}^{\alpha+2}\nonumber\\
			&=E(u_k)+\left( \frac{c_j}{c}-1\right) \left( \frac{a}{2}\|\nabla_{x,y} u_k\|_2^2+\frac{b}{4}\|\nabla_{x,y}u_k\|_2^4
			-\frac{1}{\alpha+2}\|u_k\|_{\alpha+2}^{\alpha+2}\right) \nonumber\\
			&\quad+\frac{bc_j}{4c}\bg(\frac{c_j}{c}-1\bg)\|\nabla_{x,y}u_k\|_2^4+\frac{c_j}{c}\left( 1-\bg(\frac{c_j}{c}\bg)^{\frac{\alpha}{2}}\right) \frac{1}{\alpha+2}\|u_k\|_{\alpha+2}^{\alpha+2}\nonumber\\
			&=:E(u_k)+I+II+III.\label{eq:3.2}
		\end{align}
		Using $c_j= c+o_j(1)$ and the boundedness of $u_k$ in $H_{x,y}^1$, we obtain that $I,II,III\to 0$ as $j\to\infty$. Since $\{u_k\}$ is a minimizing sequence of $m_c$, we conclude that
		\begin{align}\label{eq:3.3}
			\limsup_{j\to\infty}\,m_{c_j}\leq m_c.
		\end{align}
		 Let $u_j\in S(c_j)$ satisfy $E(u_j)\leq m_{c_j}+j^{-1}$. Using \eqref{eq:3.1} we know that $\{u_j\}$ is a bounded sequence in $H_{x,y}^1$. By the definition of $m_c$ and arguing as in \eqref{eq:3.2} we see that
		\begin{align*}
			m_c\leq E\left( \frac{\sqrt{c}}{\sqrt{c_j}}u_j\right) 
			&\leq m_{c_j}+j^{-1}+\left( \frac{c}{c_j}-1\right) \left( \frac{a}{2}\|\nabla_{x,y} u_j\|_2^2+\frac{b}{4}\|\nabla_{x,y}u_j\|_2^4
			-\frac{1}{\alpha+2}\|u_j\|_{\alpha+2}^{\alpha+2}\right) 
			\nonumber\\
			&\quad+\frac{bc}{4c_j}\bg(\frac{c}{c_j}-1\bg)\|\nabla_{x,y}u_j\|_2^4+\frac{c}{c_j}\bg(1-\bg(\frac{c}{c_j}\bg)^{\frac{\alpha}{2}}\bg)\frac{1}{\alpha+2}\|u_j\|_{\alpha+2}^{\alpha+2}.
		\end{align*}
		Taking $j\to\infty$ we then similarly deduce that
		\begin{align}\label{eq:3.4}
			m_c\leq\liminf_{j\to\infty}m_{c_j}.
		\end{align}
		Thus, we  deduce from \eqref{eq:3.3}
		and \eqref{eq:3.4} that 
		$$
			\lim_{j\to\infty}m_{c_j}=m_c,
		$$
		this implies the mapping  $c\mapsto m_c$ is continuous on $(0,\infty)$.
		
		Next, we consider the monotonicity  of $c\mapsto c^{-1} m_c$. Let $0<c_1<c_2<\infty$,  we take $u\in \mathcal{S}(c_1)$ such that $E(u)\leq m_{c_1}+\varepsilon$ for any $\varepsilon>0$. Set  $\theta\star u:=u(\theta^{-\frac{1}{d}}x,y)$, namely $\|\theta\star u\|_2^2=\theta\|u\|_2^2$, then we obtain
		\begin{align}
			m_{c_2}\leq E\left( \frac{c_2}{c_1}\star u\right) &=\frac{c_2}{c_1}E(u)+\frac{c_2}{c_1}\left[ \left( \left( \frac{c_2}{c_1}\right) ^{-\frac{2}{d}}-1\right) \|\nabla_{x,y}u\|_2^2+\left( \left( \frac{c_2}{c_1}\right) ^{1-\frac{4}{d}}-1\right) \|\nabla_{x,y}u\|_2^4\right]\nonumber\\
			& \leq \frac{c_2}{c_1}(m_{c_1}+\varepsilon).\label{eq:3.5}
		\end{align}
		By arbitrariness of $\varepsilon$, 
		\begin{align*}
			m_{c_2}\leq\frac{c_2}{c_1}m_{c_1}.
		\end{align*}
		Combine with $m_c\leq0$ for all $c\in(0,\infty)$, we infer that the mapping $c\mapsto c^{-1}m_c$ is  decreasing on $(0,\infty)$.
	\end{proof}
\end{lemma}


\begin{proof}[Proof of Theorem \ref{th:1.1}]
	\textbf{Case (i):}	 We split our proof into three steps.

\textit{\textbf{Step 1:} Existence of a  non-vanishing weak limit.}
Notice that   $\alpha<\min\{\frac{4}{d},\frac{8}{d+n}\}\leq\frac{4}{d}$, then from Proposition \ref{cor 3.1} we know that $m_c<0$ for all $c>0$. Let $\{u_k\}\subset \mathcal{S}(c)$ be a minimizing sequence of $m_c$, we then claim that there is a non-vanishing weak limit of $\{u_k\}$. By definition of $m_c$ we have
\begin{align}\label{eq:3.6}
\|u_k\|_{\alpha+2}^{\alpha+2}=(\alpha+2)\left( -m_c+o_k(1)+\frac{a}{2}\|\nabla_{x,y} u_k\|_2^2+\frac{b}{4}\|\nabla_{x,y}u_k\|_2^4\right) 
\geq (\alpha+2)(-m_c+o_k(1)).
\end{align}
Thus combining with $m_c<0$ we infer that $\liminf_{k\to\infty}\|u_k\|_{\alpha+2}>0$. The claim  follows from Lemma \ref{lem:non-vanishing}.

\textit{\textbf{Step 2:} The mapping $c\mapsto c^{-1}m_c$ is strictly decreasing on $(0,\infty)$.}
For $0<c_1<c_2<\infty$, let $\{u_k\}$ be a minimizing sequence of $m_{c_1}$. Then  we obtain
	\begin{align*}
	m_{c_2}\leq E\left( \frac{c_2}{c_1}\star u_k\right) =\frac{c_2}{c_1}E(u_k)+\frac{c_2}{c_1}\left[ \left( \left( \frac{c_2}{c_1}\right) ^{-\frac{2}{d}}-1\right) \|\nabla_{x,y}u_k\|_2^2+\left( \left( \frac{c_2}{c_1}\right) ^{1-\frac{4}{d}}-1\right) \|\nabla_{x,y}u_k\|_2^4\right],
    \end{align*}
 here $\theta\star u:=u(\theta^{-\frac{1}{d}}x,y)$. Thus combining with  $\liminf_{k\to\infty}\|u_k\|_{\alpha+2}^{\alpha+2}>0$, Lemma \ref{lem:gn1} and taking $k\to\infty$ yields
\begin{align*}
	m_{c_2}<\frac{c_2}{c_1}m_{c_1}.
\end{align*}
This completes the proof of Step 2.

\textit{\textbf{Step 3:} Conclusion.}
 Let $\{u_k\}$ be a minimizing sequence of $m_{c}$, which possesses a  $u\in H_{x,y}^1\setminus\{0\}$ such that $u_k\rightharpoonup u$ in $H_{x,y}^1$.  Using weakly lower semicontinuity of norms it is necessary that $c_1:=\|u\|_{2}^2\in(0,c]$.  Using Br\'{e}zis-Lieb lemma we have
\begin{equation*}
	E(u_k-u)+E(u)+\frac{b}{2}\|\nabla_{x,y} (u_k-u)\|_2^2\|\nabla_{x,y} u\|_2^2=E(u_k)+o_k(1),
\end{equation*}
and
\begin{equation}\label{eq:3.7}
	\|u_k-u\|_{2}^2+\|u\|_{2}^2=\|u_k\|_{2}^2+o_k(1).
\end{equation}
By \eqref{eq:3.7} we know that $\|u_k-u\|_{2}^2=c-c_1+o_k(1)$.
From  Lemma \ref{lem:3.1} we infer that
\begin{align*}
m_c=E(u)+\lim_{k\to\infty}\left( E(u_k-u)+\frac{b}{2}\|\nabla_{x,y} (u_k-u)\|_2^2\|\nabla_{x,y} u\|_2^2\right) \geq m_{c_1}+m_{c-c_1}.
\end{align*}
Therefore combining this with  Step 2 and $c_1,c-c_1<c$ we obtain
\begin{align*}
m_c>\frac{c_1}{c}m_c+\frac{c-c_1}{c}m_c=m_c,
\end{align*}
a contradiction. Thus, $\|u\|_2^2=c$ and
\begin{align*}
	m_c\leq E(u)\leq\lim_{k\to\infty}E(u_k)=m_c.
\end{align*}
So, $u$ is a minimizer of $m_c$ , and  $u$ is a critical point of $E|_{\mathcal{S}(c)}$ . 

\textbf{ Case (ii):}
In this case  we define
\begin{align*}
	c_-:=\inf\{c\in(0,\infty):m_c<0\}.
\end{align*}
For any $u \in \mathcal{S}(1)$, set
	\[
	\theta\star u:=u(\theta^{-\frac{1}{d}}x,y), \quad \theta > 0.
	\]
	Then, $\|\theta\star u\|^2_2 = \theta$, and
	\[
	E(\theta\star u) = \frac{a\theta^{1 - \frac{2}{d}}}{2} \|\nabla_{x,y} u\|_2^2 + \frac{b\theta^{2 - \frac{4}{d}}}{4} \|\nabla_{x,y} u\|_2^4 - \frac{\theta}{\alpha+2}\|u\|_{\alpha+2}^{\alpha+2}.
	\]
Since $\min\{1 - \frac{2}{d}, 2 - \frac{4}{d}\}<1$,	we have $\inf_{u\in\mathcal{S}(\theta)}E(u)= m_\theta\geq0$ for $\theta$ small enough, and then $c_->0$.
By Proposition \ref{thm:2.4},  Proposition \ref{cor 3.1} and $m_c\leq(2\pi)^{n}\widetilde{m}_{(2\pi)^{-n}c}$ we infer that $c_-\in[0,\infty)$. Hence the existence of minimizer of $m_c$ for $c\in(c_-,\infty)$ follows already from the previous proof of Theorem \ref{th:1.1} $(i)$. 

Next we  show that $m_c$ has no minimizers for $c\in(0,c_-)$. First we show that
\begin{align}\label{eq:3.8}
	m_{c_1+c_2}\leq m_{c_1}+m_{c_2},
\end{align}
for any $c_1,c_2\in(0,\infty)$. To see this, by the definition of $m_{c}$ and density arguments, for any $\vare>0$ we can find
$u_i\in C_0^\infty(\mathbb{R}^d)\otimes C^\infty_{\rm per}(\mathbb{T}^n)$  defined by \eqref{eq1.14} such that $u_i\in \mathcal{S}(c_i)$ and $m_{c_i}\leq E(u_i)+\vare$ for $i\in\{1,2\}$. Since $E(u)$ is invariant under $\mathbb{R}^d$-translation, we may also assume that $u_1,u_2$ have distinct compact support. Thus
\begin{align*}
	m_{c_1+c_2}\leq E(u_1+u_2)=E(u_1)+E(u_2)\leq m_{c_1}+m_{c_2}+2\vare.
\end{align*}
Then \eqref{eq:3.8} follows by taking $\vare\to0$. Combining with $m_c\leq0$,  we know that $c\mapsto c^{-1}m_c$ is monotone decreasing.  From Lemma \ref{lem:3.1},   the mapping $c\mapsto m_c$ is continuity on $(0,\infty)$, then $m_c=0$ for all $c\in(0,c_-]$. Now let $c_3\in(0,c_-)$, we suppose that  $m_{c_3}$ has an minimizer $u$, then  arguing as \eqref{eq:3.5}, we have
\begin{align*}
	m_\tau\leq E\left( \frac{\tau}{c_3}\star u\right) &=\bg(\frac{\tau}{c_3}\bg)\left( \bg(\frac{\tau}{c_3}\bg)^{-\frac{2}{d}}\frac{a}{2}\|\nabla_{x,y}u\|_2^2
	+\bg(\frac{\tau}{c_3}\bg)^{1-\frac{4}{d}}\frac{b}{4}\|\nabla_{x,y}u\|_2^4
	-\frac{1}{\alpha+2}\|u\|_{\alpha+2}^{\alpha+2}\right) \\
	&<\bg(\frac{\tau}{c_3}\bg)E(u)= \bg(\frac{\tau}{c_3}\bg)m_{c_3}=0.
\end{align*}
This in turn implies that $m_\tau< m_{c_3}$, which also contradicts the fact that $m_c=0$ for all $c\in(0,c_-]$ that we deduced previously.

It remains to show that $m_{c_-}$ is achieved by some $v \in\mathcal{S}(c_-)$ when $\alpha\in(\frac{4}{d},\frac{8}{d+n})$. Let $c_k := c_- + 1/k$ for any $k \in \mathbb{N}^+$. Since $m_{c_k} < 0$ by the define of $c_-$,  one may choose $v_k \in \mathcal{S}(c_k)$ such that
\begin{equation}\label{eq:v_k}
	 E(v_k) =m_{c_k} < 0 \qquad \text{for each}~ k \in \mathbb{N}^+.
\end{equation}
Hence the sequence $\{v_k\}$ is bounded in $H_{x,y}^1$ via \eqref{eq:3.1}, and  $\lim_{k \to \infty}\|\nabla_{x,y} v_k\|_2^2 $ and  $\lim_{k \to \infty} \|v_k\|_{\alpha+2}^{\alpha+2}$ exist (up to a subsequence). In view of \eqref{eq:v_k}, it is also clear that $\lim_{k \to \infty}\|\nabla_{x,y} v_k\|_2^2 >0$. 
Indeed, if $\|\nabla_{x,y} v_k\|_2^2\to0$, from Lemma \ref{lem:gn2} (since $\alpha\in(\frac{4}{d},\frac{8}{d+n})$ and $d+n\leq4$, we have $n=1$) we infer that
\begin{align*}
	E(v_k) &\geq  \frac{a}{2}\|\nabla_{x,y} v_k\|_2^2 + \frac{b}{4} \|\nabla_{x,y} v_k\|_2^4 - \frac{C(\alpha,d)}{\alpha+2}\|\nabla_x v_k\|_2^{\frac{\alpha d}{2}}\|v_k\|_2^{\frac{4-\alpha(d-1)}{2}}
	\left( \| v_k\|_{2}^{\frac{\alpha}{2}}+\|\pt_y v_k\|_{2}^{\frac{\alpha}{2}}\right)\\
	&\geq  \frac{a}{2}\|\nabla_{x,y} v_k\|_2^2 + \frac{b}{4} \|\nabla_{x,y} v_k\|_2^4 - C(\alpha,d,c)\|\nabla_{x,y} v_k\|_2^{\frac{\alpha d}{2}}\geq0.
\end{align*}
This contradicts \eqref{eq:v_k}.
We claim that $\{v_k\}$ is non-vanishing.
Otherwise, $\|v_k\|_{\alpha+2}^{\alpha+2}\to0$. Then we obtain
\begin{equation*}
	\lim_{k \to \infty} E(v_k) \geq \frac{a}{2}\lim_{k \to \infty}\|\nabla_{x,y} v_k\|^2_{2} > 0.
\end{equation*}
This contradicts \eqref{eq:v_k} and hence $\{v_k\}$ is non-vanishing. 	Then by Lemma \ref{lem:non-vanishing} there exist $\{x_j\}\subset \R^d$ and some $v\in H_{x,y}^1\setminus\{0\}$ such that up to a subsequence
$	v_k(x+x_j,y)\rightharpoonup v(x,y)$ in $H_{x,y}^1$. Denote $c' := \|v\|^2_{2} \in (0, c_-]$ and $w_k := v_k(\cdot + x_j) - v$. Clearly, $\lim_{k \to \infty} \|w_k\|^2_{2} = c_- -c'$ and
\begin{equation*}
	\lim_{k \to \infty} \|\nabla_{x,y} v_k\|^2_{2} = \|\nabla_{x,y} v\|^2_{2} + \lim_{k \to \infty} \|\nabla_{x,y} w_k\|^2_{2}.
\end{equation*}
Noting that $m_{c_k} \to m_{c_-} = 0$, we deduce from \eqref{eq:v_k} that
\begin{equation}\label{eq:Iv_k}
	0 = \lim_{k \to \infty} E(v_k) = E(v) + \lim_{k \to \infty}E(w_k)+\lim_{k \to \infty}\frac{b}{2}\|\nabla_{x,y} (v_k-v)\|_2^2\|\nabla_{x,y} v\|_2^2.
\end{equation}
Since $m_s = 0$ for any $s \in (0, c_-]$,   $E(v) \geq 0$ and $\lim_{k \to \infty} E(w_k) \geq 0$. In view also of \eqref{eq:Iv_k}, we obtain
\begin{equation*}
E(v) = 0 \qquad \text{and} \qquad \lim_{k \to \infty} E(w_k) = 0.
\end{equation*}
Now, from the fact that   $E_s = 0$ is not  achieved  for any $s \in (0, c_-)$, it follows that $c' = c_-$ and thus $v \in \mathcal{S}(c_-)$ is a minimizer of $m_{c_-} = 0$. As a consequence, $v \in \mathcal{S}(c_-)$ is an energy ground state solution  of problem \eqref{eq:1.1} when $c = c_-$.
\end{proof}

\subsection{$y$-dependence of normalized ground states}\label{sub1.2}

The lack of scaling invariance in problem \eqref{eq:1.1} poses a key challenge. To prove Theorem \ref{th:1.2}, we adapt the method of \cite{Luo2022arXiv,Tzvetkov2014AnalPDE} to verify \eqref{tau lambda}.
\begin{lemma}\label{lem:4.1}
	Let $\alpha\in(2,\min\{\frac{8}{d+n},\frac{4}{d}\})$, then we have $\lim_{\ld\to0}m_{1,\ld}<(2\pi)^n\widehat{m}_{(2\pi)^{-n}}$, here $m_{c,\ld}$ is defined in \eqref{m minimizing lambda}.
\end{lemma}

\begin{proof} 
	Let $\delta\in(0,\pi)$ and $\beta\in(0,\infty)$, we define the function $\rho:[0,2\pi]\to[0,\infty)$ :
	\begin{align}\label{eq rho}
		\rho(z)=\left\{
		\begin{array}{ll}
			0,&z\in[0,\delta]\cup[2\pi-\delta,2\pi],\\
			\beta(z-\delta),&z\in[\delta,\pi],\\
			\beta(2\pi-\delta-z),&z\in[\pi,2\pi-\delta].
		\end{array}
		\right.
	\end{align}
	Direct calculation shows that
	\begin{align*}
		\|\rho\|_{L^{\alpha+2}(\mathbb{T})}^{\alpha+2}=\dfrac{2\beta^{\alpha+2}(\pi-\delta)^{\alpha+3}}{\alpha+3}\quad\text{and}\quad\|\rho\|_{L^2(\mathbb{T})}^2=\dfrac{2\beta^2(\pi-\delta)^{3}}{3}.
	\end{align*}
	Then, we chose the value of $\beta$ by
	$$
	\beta=2(\pi-\delta)^{-1}\left( \frac{\alpha+3}{3}\right) ^{\frac{1}{\alpha}},
	$$ 
	such that 
	$$\|\rho\|_{L^{\alpha+2}(\mathbb{T})}^{\alpha+2}>\|\rho\|_{L^2(\mathbb{T})}^2= \frac{8}{3}\bg(\frac{\alpha+3}{3}\bg)^{\frac{2}{\alpha}}(\pi-\delta). $$
	Using Proposition \ref{thm:2.4} and \eqref{invariant hat minimizing} below, we know that $\widehat{m}_{\|\rho\|^{-2n}_{L^2(\mathbb{T})},\lambda}<0$  has an minimizer $Q_\ld$ for all $0<\ld\ll 1$ with $\lambda=c^{\frac{2(\alpha-2)}{\alpha d-8}}$. Since  the mapping $\ld\mapsto \widehat{m}_{c,\lambda}$ is monotone increasing for fixed $c>0$, there exists some $A>0$ and $\ld_0\in(0,\infty)$ such that $\widehat{m}_{\|\rho\|^{-2n}_{L^2(\mathbb{T})},\lambda}\leq -A$ for all $\ld\in(0,\ld_0)$. Arguing as in \eqref{eq:3.6} we then infer that there exists some $C_A>0$ such that
	\begin{align}\label{eq4.10}
		\inf_{\ld\in(0,\ld_0)}\|Q_\ld\|_{L_x^{\alpha+2}}^{\alpha+2}\geq C_A>0.
	\end{align}
	Next, we show that $m_{1,\ld}<(2\pi)^{n} \widehat{m}_{(2\pi)^{-n},\lambda}$ for all $\ld\ll 1$. 
	Let $\rho_\vare$ be the $\vare$-mollifier of $\rho$ on $[0,2\pi]$ with some to be determined small $\vare>0$. In particular,  $\rho$ has compact support in $(0,2\pi)$, so is $\rho_\vare$ for $\vare\ll 1$. Next, we define the function $\psi$ by
	\begin{align*}
		\psi(x,y)&:=Q_\ld(x)\left( \|\rho\|_{L^2(\mathbb{T})}/\|\rho_\vare\|_{L^2(\mathbb{T})}\right) ^{n}\Psi_\vare(y),
	\end{align*}
	where $\Psi_\vare(y)=\prod_{j=1}^n\rho_\vare(y_j)$. This is to be understood that we extend $\rho_\vare$ $2\pi$-periodically along the $y$-direction, which is possible since $\rho_\vare$ has compact support in $(0,2\pi)$ when $\vare\ll 1$. We have then $\psi\in \mathcal{S}(1)$.	Consequently,
	\begin{align}
		m_{1,\ld}\leq E_\ld(\psi)
		=&\|\rho\|_{L^2(\mathbb{T})}^{2n}\left( \frac{a\lambda^{\frac{4-d}{\alpha-2}+1-\frac{d}{2}}}{2}\|\nabla_x Q_\ld\|_{L_x^2}^2+\frac{b}{4}\|\nabla_x Q_\ld\|_{L_x^2}^4
		-\frac{1}{\alpha+2}\|Q_\ld\|_{L_x^{\alpha+2}}^{\alpha+2}\right) \nonumber\\
		&+\frac{b}{4}\|\rho\|_{L^2(\mathbb{T})}^{2n}\left( \|\rho\|_{L^2(\mathbb{T})}^{2n}-1\right) \|\nabla_x Q_\ld\|_{L_x^2}^4\nonumber\\
		&-\frac{\|\rho\|_{L^2(\mathbb{T})}^{2n}}{\alpha+2}
		\left( \|\rho\|_{L^2(\mathbb{T})}^{-2n}\|\rho\|_{L^2(\mathbb{T})}^{(\alpha+2)n}\|\rho_\vare\|_{L^{2}(\mathbb{T})}^{-(\alpha+2)n}\|\rho_\vare\|_{L^{\alpha+2}(\mathbb{T})}^{(\alpha+2)n}-1\right) 
		\|Q_\ld\|_{L_x^{\alpha+2}}^{\alpha+2}\nonumber\\
		&+ \frac{b\lambda^2}{4\|\rho_\vare\|_{L^2(\mathbb{T})}^{4n}}\|\nabla_y\Psi_\vare\|_{L_y^2}^4+\frac{b\lambda\|\rho\|_{L^2(\mathbb{T})}^{2n}}{2\|\rho_\vare\|_{L^2(\mathbb{T})}^{2n}}\|\nabla_{x} Q_\lambda\|^2_{L_x^2}\|\nabla_{y} \Psi_\varepsilon\|^2_{L_y^2}
		+\frac{a\lambda^{\frac{4-d}{\alpha-2}+2-\frac{d}{2}}}{2\|\rho_\vare\|_{L^2(\mathbb{T})}^{2n}}\|\nabla_y \Psi_\vare\|_{L_y^2}^2 \nonumber\\
		=&\|\rho\|_{L^2(\mathbb{T})}^{2n}\widehat{m}_{\|\rho\|_{L^2(\mathbb{T})}^{-2n},\lambda}+I_1+I_2+I_3.\nonumber
	\end{align}
	As in Lemma \ref{lem:3.1}, there exist $c_{0,\lambda}$ such that the function $c\mapsto c^{-1}\widehat{m}_{c,\lambda}$ is strictly decreasing on $(c_{0,\lambda},\infty)$. 	Let $\delta\in(0,\pi)$ such that $\delta>\pi-\frac{3}{8}\bg(\frac{3}{\alpha+3}\bg)^{\frac2\alpha}$, by direct calculation we see that $\rho\in H^1(\mathbb{T})$ and $1>\|\rho\|_{L^2(\mathbb{T})}^2$. Then $\|\rho\|_{L_y^2}^{2n}\widehat{m}_{\|\rho\|_{L_y^2}^{-2n},\lambda}
	<(2\pi)^n\widehat{m}_{(2\pi)^{-n},\lambda}$ and $I_1\leq0$.  By the properties of a mollifier operator we also know that
	\begin{align*}
		\|\rho_\vare\|^{2}_{L^2(\mathbb{T})}=\|\rho\|_{L^2(\mathbb{T})}^2+o_\vare(1)\quad\text{and}\quad\|\rho_\vare\|_{L^{\alpha+2}(\mathbb{T})}^{\alpha+2}
		=\|\rho\|_{L^{\alpha+2}(\mathbb{T})}^{\alpha+2}+o_\vare(1).
	\end{align*}
	Note  $\|\rho\|_{L^2(\mathbb{T})}^2<\|\rho\|_{L^{\alpha+2}(\mathbb{T})}^{\alpha+2}$, we have
	\begin{align}
		\|\rho\|_{L^2(\mathbb{T})}^{-2n}\|\rho\|_{L^2(\mathbb{T})}^{(\alpha+2)n}\|\rho_\vare\|_{L^{2}(\mathbb{T})}^{-(\alpha+2)n}\|\rho_\vare\|_{L^{\alpha+2}(\mathbb{T})}^{(\alpha+2)n}-1
		=\frac{\|\rho\|_{L^2(\mathbb{T})}^{(\alpha+2)n}(\|\rho\|^{(\alpha+2)n}_{L^{\alpha+2}(\mathbb{T})}+o_\vare(1))}
		{\|\rho\|_{L^2(\mathbb{T})}^{2n}(\|\rho\|^{(\alpha+2)n}_{L^2(\mathbb{T})}+o_\vare(1))}-1 \geq \zeta>0\label{eq4.11}
	\end{align}
	for some $\zeta>0$ and all $\vare\ll 1$. We then fix some $\vare$ so that \eqref{eq4.11} holds and combining with \eqref{eq4.10}
	\begin{align*}
		I_2\leq -\frac{\|\rho\|_{L^2(\mathbb{T})}^{2n}C_A\zeta}{\alpha+2}.
	\end{align*}
	Consequently, for fixed $\delta$ and $\vare$ we see that $\lim_{\ld\to 0}I_3=0 $. The desired claim follows by taking $\ld$ small.
\end{proof}
Proceeding as in the proof of Theorem \ref{th:1.1} we know that  $\mu_{1,\tau}$ (defined in \eqref{mu c tau}) has an minimizer $u_\tau\in \mathcal{S}(1)$ for any $\tau\in(0,\infty)$.
\begin{lemma}\label{lem:4.2}
	Let  $\alpha\in(0,\frac{4}{d+n})$ and $u_\tau$ be a minimizer of $\mu_{1,\tau}$. Then we have
	\begin{align}
		\lim_{\tau\to\infty}\mu_{1,\tau}=(2\pi)^n\widehat{\mu}_{(2\pi)^{-n},\infty}\label{eq:4.4},
	\end{align}
	and
	\begin{align}
		&\lim_{\tau \to \infty}\tau^2\|\nabla_{y} u_{\tau}\|^2_{2}=0,\label{eq:4.5} \\
		&\lim_{\tau \to \infty}\tau^{2+d-\frac{4}{\alpha}}\|\nabla_{y} u_{\tau}\|^4_{2}=0,\\
		&\lim_{\tau \to \infty}\frac{1}{\tau^{\frac{4}{\alpha}-d}}\|\nabla_{x} u_{\tau}\|^2_{2}\|\nabla_{y} u_{\tau}\|^2_{2}=0 ,
		\\
		&\lim_{\tau \to \infty}\frac{1}{\tau^{\frac{4}{\alpha}+2-d}}\|\nabla_{x} u_{\tau}\|^4_{2}=0,\label{eq:4.6}
	\end{align}
	where $\mu_{c,\tau}$ and $\widehat{\mu}_{c,\infty}$ are defined in \eqref{mu c tau} and \eqref{h m c} respectively.
\end{lemma}
\begin{proof}
	Let $u\in H^1_x$ be such that $\|u\|_{L_x^2}^2 = (2\pi)^{-n}$ and $\widehat{E}^{\infty}(u) = \widehat{\mu}_{(2\pi)^{-n},\infty}$. Then $u \in H^1_{x,y}$, $\|u\|_{2}^2 = 1$ and
	\begin{align}
		\mu_{1,\tau} \leq E^{\tau}(u) &= (2\pi)^n\left(\frac{a}{2}\|\nabla_{x} u\|^2_{L_x^2}+\frac{b(2\pi)^n}{4\tau^{\frac{4}{\alpha}+2-d}}\|\nabla_{x} u\|^4_{L_x^2} -\frac{1}{\alpha+2}\|u\|^{\alpha+2}_{L_x^{\alpha+2}}\right)\nonumber\\
		& =(2\pi)^n\left( \widehat{E}^{\infty}(u)+\frac{b(2\pi)^n}{4\tau^{\frac{4}{\alpha}+2-d}}\|\nabla_{x} u\|^4_{L_x^2}\right) = (2\pi)^n\widehat{\mu}_{(2\pi)^{-n},\infty} + o_\tau(1),\label{eq:4.7}
	\end{align}
	where $\lim_{\tau \to \infty}o_\tau(1) = 0$, implying that 
	\begin{equation}\label{eq:4.8}
		\limsup_{\tau \to \infty}\mu_{1,\tau} \leq (2\pi)^n\widehat{\mu}_{(2\pi)^{-n},\infty}.
	\end{equation}
	Note that
	\begin{align}
		\mu_{1,\tau} = E^{\tau}(u_{\tau}) =&\frac{a(\tau^2-1)}{2}\|\nabla_{y} u_{\tau}\|^2_{2}+\frac{b}{2\tau^{\frac{4}{\alpha}-d}}\|\nabla_{x} u_{\tau}\|^2_{2}\|\nabla_{y} u_{\tau}\|^2_{2}+\frac{b\tau^{2+d-\frac{4}{\alpha}}}{4}\|\nabla_{y} u\|^4_{2}
		+\frac{b}{4\tau^{\frac{4}{\alpha}+2-d}}\|\nabla_{x} u\|^4_{2}\nonumber\\
		&+\frac{a}{2}\|\nabla_{y} u_{\tau}\|^2_{2}+\frac{a}{2}\|\nabla_{x} u_{\tau}\|^2_{2}
		-\frac{1}{\alpha+2}\|u_{\tau}\|^{\alpha+2}_{\alpha+2}.\nonumber
	\end{align}
	Since $\alpha \in (0,\frac{4}{d+n})$,  by the  Gagliardo-Nirenberg inequality (Lemma \ref{lem:gn1}), we obtain
	\begin{align*}
		\inf_{\tau>0}\left(\frac{a}{2} +\frac{a}{2}\|\nabla_{x,y} u_{\tau}\|^2_{2} -\frac{1}{\alpha+2}\|u_{\tau}\|^{\alpha+2}_{\alpha+2} \right)&\geq\inf_{\tau>0}\left(\frac{a}{2}\| u_{\tau}\|^2_{H^1_{x,y}} -C\|u_{\tau}\|^{\frac{\alpha(d+n)}{2}}_{H^1_{x,y}} \right) \\
		&\geq\inf_{t>0}\left( C_1t^2-C_2 t^{\frac{\alpha(d+n)}{2}}\right) >-\infty.
	\end{align*}
	Hence,
	\begin{align}
		&0 \leq \limsup_{\tau \to \infty}(\tau^2-1)\|\nabla_{y} u_{\tau}\|^2_{2} <\infty, \nonumber\\
		&0 \leq \limsup_{\tau \to \infty}\frac{1}{\tau^{\frac{4}{\alpha}-d}}\|\nabla_{x} u_{\tau}\|^2_{2}\|\nabla_{y} u_{\tau}\|^2_{2} <\infty, \nonumber\\
		&0 \leq \limsup_{\tau \to \infty}\tau^{2+d-\frac{4}{\alpha}}\|\nabla_{y} u_\tau\|^4_{2}<\infty,\nonumber\\
		&0 \leq\limsup_{\tau \to \infty}\frac{1}{\tau^{\frac{4}{\alpha}+2-d}}\|\nabla_{x} u_{\tau}\|^4_{2} <\infty\nonumber,
	\end{align}
	showing that
	\begin{align}\label{eq:4.9}
		\lim_{\tau\to\infty}\|\nabla_{y} u_{\tau}\|^2_{2}=0.
	\end{align}
	Observe that estimating similarly as in \eqref{eq:3.1}, if $\|u_\tau\|_{H_{x,y}^1}\to \infty$ as $\tau\to\infty$ would in turn imply the contradiction $\mu_{1,\tau}\to\infty$, thus $\{u_\tau\}$ is uniformly bounded in $H_{x,y}^1$
	with respect to $\tau$.
 Defining
	\begin{align*}
		w_\tau(y):=\|u_\tau(\cdot,y)\|_{L_x^2}^2.
	\end{align*}
	Using the product rule, H\"older inequality and \eqref{eq:4.9} we also infer that
	\begin{equation}\label{eq:4.10}
		\| \nabla_yw_\tau \|_{L^1_{y}}\leq C \| u_\tau \|_2\| \nabla_yu_\tau\|_2 \to 0\quad\text{as} \quad\tau\to\infty,
	\end{equation}
	we then conclude that $\{w_\tau\}$ is a bounded sequence in $W^{1,1}(\mathbb{T}^n)$ with $\|w_\tau\|_{L_y^1}=1$. By a rescaling arguments it is easy to show that
	\begin{align}\label{invariant hat mu}
	\widehat{\mu}_{c,\infty}=c^{\frac{d-2-4/\alpha}{d-4/\alpha}}\widehat{\mu}_{1,\infty}.
	\end{align} 
	Therefore, from \eqref{invariant hat mu} and \eqref{eq:4.7}
	\begin{align}
		&\mu_{1,\tau}=E^{\tau}(u_\tau)\geq\frac{a}{2}\|\nabla_{x} u_\tau\|^2_{2}
		-\frac{1}{\alpha+2}\|u_\tau\|^{\alpha+2}_{\alpha+2}= \int_{y\in\mathbb{T}^n}\widehat{E}^{\infty}(u_\tau(\cdot,y))\,dy\nonumber
		\geq\int_{y\in\mathbb{T}^n}\widehat{\mu}_{\|u_\tau(\cdot,y)\|^2_{L_x^2},\infty}\,dy\nonumber\\&=\int_{y\in\mathbb{T}^n}(\|u_\tau(\cdot,y)\|^2_{L_x^2})^{\frac{d-2-4/\alpha}{d-4/\alpha}}\widehat{\mu}_{1,\infty}\,dy
		=\|w_\tau\|_{L_y^{\frac{d-2-4/\alpha}{d-4/\alpha}}}^{\frac{d-2-4/\alpha}{d-4/\alpha}}\widehat{\mu}_{1,\infty}.\label{eq:4.11}
	\end{align}
	Using $\alpha\in (0,\frac{4}{d+n})$ we deduce that $\frac{d-2-4/\alpha}{d-4/\alpha}\in(1,1+\frac{2}{n})$.  When $n\leq 2$, by Rellich's compact embedding we have $W_y^{1,1}\hookrightarrow \hookrightarrow  L_y^{r}$ is compact embedding for all $r\in[1,1+\frac{2}{n})$. Using \eqref{eq:4.10} and $\|w_\tau\|_{L_y^1}= 1$ we deduce
	\begin{align}
		\lim_{\tau\to\infty}\|w_\tau\|_{L_y^{\frac{d-2-4/\alpha}{d-4/\alpha}}}^{\frac{d-2-4/\alpha}{d-4/\alpha}}
		=(2\pi)^n((2\pi)^n)^{\frac{2+4/\alpha-d}{d-4/\alpha}}=(2\pi)^{\frac{2n}{d-4/\alpha}}.\label{eq:4.12}
	\end{align}
	For $n\geq 3$, we use Minkowski's inequality and the Sobolev embedding $H^1_y\hookrightarrow L_y^{\frac{2n}{n-2}}$ to obtain
	\begin{align*}
		\sup_{\tau>0}\|w_\tau\|_{L_y^{\frac{n}{n-2}}}= \sup_{\tau>0}\|u_\tau\|^2_{L_y^{\frac{2n}{n-2}}L_x^2}
		\leq C\sup_{\tau>0}\|u_\tau\|^2_{L_x^2L_y^{\frac{2n}{n-2}}}\leq C
		\sup_{\tau>0}\|u_\tau\|^2_{L_x^2H_y^1}<\infty.
	\end{align*}
	Combining with $W_y^{1,1}\hookrightarrow \hookrightarrow L_y^{1}$ and interpolation we arrive at \eqref{eq:4.12} again.
	Thus using \eqref{invariant hat mu}, \eqref{eq:4.11} and \eqref{eq:4.12}, we infer that
	\begin{align}
		\liminf_{\tau\to\infty} \mu_{1,\tau}\geq (2\pi)^n\left( ((2\pi)^n)^{\frac{2+4/\alpha-d}{d-4/\alpha}}\widehat{\mu}_{1,\infty}\right) =(2\pi)^n\widehat{\mu}_{(2\pi)^{-n},\infty}.\label{eq:4.13}
	\end{align}
	Then \eqref{eq:4.4} follows from \eqref{eq:4.8} and \eqref{eq:4.13}. Finally, by keeping 
	$$
	\frac{a\tau^2}{2}\|\nabla_{y} u_\tau\|^2_{2}
	+\frac{b}{2\tau^{\frac{4}{\alpha}-d}}\|\nabla_{x} u_{\tau}\|^2_{2}\|\nabla_{y} u_{\tau}\|^2_{2}+\frac{b\tau^{2+d-\frac{4}{\alpha}}}{4}\|\nabla_{y} u_\tau\|^4_{2}
	+\frac{b}{4\tau^{\frac{4}{\alpha}+2-d}}\|\nabla_{x} u_\tau\|^4_{2}
	$$ 
	in \eqref{eq:4.11} we have
	\begin{align}
		\mu_{1,\tau}\geq& \frac{a\tau^2}{2}\|\nabla_{y} u_\tau\|^2_{2}
		+\frac{b}{2\tau^{\frac{4}{\alpha}-d}}\|\nabla_{x} u_{\tau}\|^2_{2}\|\nabla_{y} u_{\tau}\|^2_{2}+\frac{b\tau^{2+d-\frac{4}{\alpha}}}{4}\|\nabla_{y} u_\tau\|^4_{2}\nonumber\\
		&+\frac{b}{4\tau^{\frac{4}{\alpha}+2-d}}\|\nabla_{x} u_\tau\|^4_{2}+\|w_\tau\|_{L_y^{\frac{d-4/\alpha-2}{d-4/\alpha}}}^{\frac{d-4/\alpha-2}{d-4/\alpha}}\widehat{\mu}_{1,\infty}.\label{eq:4.14}
	\end{align}
	Then, \eqref{eq:4.5}-\eqref{eq:4.6}  follow  from \eqref{eq:4.4}, \eqref{eq:4.12} and taking $\tau$ in \eqref{eq:4.14} to infinity.
\end{proof}

\begin{lemma}\label{lem:4.3}
	Let  $\alpha\in(0,\frac{4}{d+n})$ and $u_\tau$ be a minimizer of $\mu_{1,\tau}$. For any $\tau>0$, then we have
	\begin{align}
		a\|\nabla_x u_\tau\|_2^2+\frac{b}{\tau^{\frac{4}{\alpha}-d}}\|\nabla_{x} u_{\tau}\|^2_{2}\|\nabla_{y} u_{\tau}\|^2_{2}+\frac{b}{\tau^{\frac{4}{\alpha}+2-d}}\|\nabla_x u_\tau\|_2^4 -\frac{\alpha d}{2(\alpha+2)}\|u_\tau\|_{\alpha+2}^{\alpha+2}=0.\label{eq:4.15}
	\end{align}
	Consequently, there exists some $\omega_\tau\in\R$ such that
	\begin{align}
		&-\left( a+\frac{b}{\tau^{\frac{4}{\alpha}+2-d}} \|\nabla_{x}u_\tau\|_2^2+\frac{b}{2\tau^{\frac{4}{\alpha}-d}}\|\nabla_{y}u_\tau\|^2_2\right) \Delta_{x}u_\tau\nonumber\\
		&-\left( a\tau^2+\frac{b}{2\tau^{\frac{4}{\alpha}-d}}\|\nabla_{x}u_\tau\|^2_2+b\tau^{2+d-\frac{4}{\alpha}}\|\nabla_{y}u_\tau\|^2_2\right) \Delta_{y}u_\tau
		+\omega_\tau u_\tau=|u_\tau|^\alpha u_\tau.\label{eq:4.16}
	\end{align}
	Moreover, we have 
	\begin{align*}
		\lim_{\tau\to\infty}\omega_\tau=\hat{\omega}:=-\frac{2(\alpha d-2\alpha-4)}{\alpha d-4} (2\pi)^n \widehat{\mu}_{(2\pi)^{-n},\infty} .
	\end{align*}
\end{lemma}

\begin{proof}
	Since $E^\tau(u_\tau)=\mu_{1,\tau}$, we have $\frac{d}{dt}E^\tau(u_\tau\ast t)|_{t=1}=0$, here 
	the scaling operator $u\mapsto u\ast t$ for $t\in(0,\infty)$ is defined by
	\begin{align}\label{eq:L2scaling}
		(u\ast t)(x,y):=t^{\frac{d}{2}}u(tx,y)\quad\text{for all $u\in H^1_{x,y}$}.
	\end{align} 
	Then
	\eqref{eq:4.15} follows from $\frac{d}{dt}E^\tau(u_\tau\ast t)|_{t=1}=0$, and \eqref{eq:4.16} is a direct consequence of Lagrange multiplier theorem. Next, testing \eqref{eq:4.16} with $u_\tau$ and combining with \eqref{eq:4.15}, we get
	\begin{align*}
		\omega_\tau=&-a\left( 1-\frac{2(\alpha+2)}{d\alpha}\right) \|\nabla_x u_\tau\|_2^2-a\tau^2\|\nabla_y u_\tau\|_2^2-\frac{b}{\tau^{\frac{4}{\alpha}+2-d}}\left( 1-\frac{2(\alpha+2)}{d\alpha}\right) \|\nabla_x u_\tau\|_2^4\\
		&-\frac{b}{\tau^{\frac{4}{\alpha}-d}}\left( 1-\frac{2(\alpha+2)}{d\alpha}\right)\|\nabla_{x} u_{\tau}\|^2_{2}\|\nabla_{y} u_{\tau}\|^2_{2}-b\tau^{2+d-\frac{4}{\alpha}}\|\nabla_y u_\tau\|_2^4,
	\end{align*}
	by \eqref{eq:4.5}-\eqref{eq:4.6} we obtain
	\begin{equation}\label{eq:4.18}
		\omega_\tau=-a\left( 1-\frac{2(\alpha+2)}{d\alpha}\right) \|\nabla_x u_\tau\|_2^2+o_\tau(1),
	\end{equation}
	where $\lim_{\tau\to\infty}o_\tau(1)=0.$
	On the other hand,  by \eqref{eq:4.15}, we get
	\begin{align*}
		\mu_{1,\tau}=E^\tau(u_\tau)=&\frac{a}{2}\left( 1-\frac{4}{d\alpha}\right) \|\nabla_x u_\tau\|_2^2+\frac{a\tau^2}{2}\|\nabla_y u_\tau\|_2^2+\frac{b}{4\tau^{\frac{4}{\alpha}+2-d}}\left( 1-\frac{8}{d\alpha}\right) \|\nabla_x u_\tau\|_2^4\\
		&+\frac{b}{2\tau^{\frac{4}{\alpha}-d}}\left( 1-\frac{4}{d\alpha}\right)\|\nabla_{x} u_{\tau}\|^2_{2}\|\nabla_{y} u_{\tau}\|^2_{2}+\frac{b\tau^{2+d-\frac{4}{\alpha}}}{4}\|\nabla_y u_\tau\|_2^4,
	\end{align*}
	then by \eqref{eq:4.5}-\eqref{eq:4.6}, we infer that
	\begin{equation}\label{eq:4.19}
		\mu_{1,\tau}=\frac{a}{2}\left( 1-\frac{4}{d\alpha}\right) \|\nabla_x u_\tau\|_2^2+o_\tau(1).
	\end{equation}
	Combining with \eqref{eq:4.4}, \eqref{eq:4.18} and \eqref{eq:4.19} we obtain
	\begin{align*}
		\omega_\tau=-\frac{2(\alpha d-2\alpha-4)}{\alpha d-4}(2\pi)^n \widehat{\mu}_{(2\pi)^{-n},\infty}+o_\tau(1),
	\end{align*}
	and the desired claim follows.
\end{proof}

\begin{proposition}\label{lem:4.4}
	Let  $\alpha\in(0,\frac{4}{d+n})$. There exists a $v \in \widehat{\mathcal{S}}((2\pi)^{-n})$ and $\{x_k\}\subset\mathbb{R}^d$ such that, up to subsequence and $\mathbb{R}^d$-translations, $v_\tau:=u_\tau(x+x_k,y)\to v$  in $H_{x,y}^1$  as $\tau\to\infty$, $\widehat{E}^{\infty}(v)=\widehat{\mu}_{(2\pi)^{-n},\infty}$ and $v$ satisfies
	\begin{align}
		-a \Delta_{x}v +\hat{\omega} v=|v|^\alpha v.\label{eq:4.20}
	\end{align}
\end{proposition}
\begin{proof}
	From the proof of Lemma \ref{lem:4.2} we know that $\{u_\tau\}$ is a bounded sequence in $H_{x,y}^1$. Moreover, using \eqref{eq:4.4}-\eqref{eq:4.6} and \eqref{eq:4.15} we obtain
	\begin{align*}
		\lim_{\tau \to \infty}\|u_\tau\|_{\alpha+2}^{\alpha+2}\geq\frac{4(\alpha+2)}{\alpha d-4}(2\pi)^n\widehat{\mu}_{(2\pi)^{-n},\infty}>0.
	\end{align*}
	Then by Lemma \ref{lem:non-vanishing}, there exists $\{x_k\}\subset\mathbb{R}^d$ such that, up to subsequence and $\mathbb{R}^d$-translations, $v_\tau:=u_\tau(x+x_k,y)\rightharpoonup v\in H_{x,y}^1\setminus\{0\}$ in $H_{x,y}^1$.  Moreover,  $v_\tau\rightharpoonup v$ in $H_{x,y}^1$, the Sobolev embedding, \eqref{eq:4.5}-\eqref{eq:4.6}, \eqref{eq:4.16} and $\lim_{\tau\to\infty}\omega_\tau=\hat{\omega}$ yield \eqref{eq:4.20}.
	
	Next, we show that $v$ is a minimizer of $\widehat{\mu}_{(2\pi)^{-n},\infty}$.	 From \eqref{eq:4.5} we get $\nabla_yv = 0$. In particular $v$ is $y$-independent.
	By weakly lower 	 semicontinuity of norm  we have $\|v\|_{L_x^2}^2\in(0,(2\pi)^{-n}]$. We claim that $\|v\|_{L_x^2}^2=(2\pi)^{-n}$. Assume that $\|v\|_{L_x^2}^2:=c_1\in (0,(2\pi)^{-n})$. Using Br\'{e}zis-Lieb lemma, \eqref{eq:4.4}, \eqref{eq:4.5}-\eqref{eq:4.6} we have
	\begin{align}
		&E^\infty(v_\tau-v)+E^\infty(v)=(2\pi)^n\widehat{\mu}_{(2\pi)^{-n},\infty}+o_\tau(1),\nonumber\\
		&\|v_\tau-v\|_{2}^2+(2\pi)^nc_1=1+o_\tau(1).\label{eq:4.21}
	\end{align}
	Define $w_\tau:=\|v_\tau-v\|_{L_{x}^2}^2$, using \eqref{eq:4.10} and \eqref{eq:4.21} we have
	\begin{align*}
		\|w_\tau\|_{L_y^1}=1-(2\pi)^nc_1+o_\tau(1)~~~~\text{ and }~~~~	\|\nabla_y w_\tau\|_{L_y^1}=o_\tau(1).
	\end{align*}
	Then arguing as in \eqref{eq:4.13} we infer that
	$$ 
	\liminf_{\tau\to\infty} E^\tau(v_\tau-v)\geq (2\pi)^n\widehat{\mu}_{(2\pi)^{-n}-c_1,\infty}.
	$$
	By definition we also know that $E^\infty(v)\geq (2\pi)^n\widehat{\mu}_{c_1,\infty}$. Arguing as Lemma \ref{lem:3.1} and the proof of Theorem \ref{th:1.1} Step 2, we deduce that the mapping $c\mapsto c^{-1}\widehat{\mu}_{c,\infty}$ is strictly decreasing on $(0,\infty)$. Summing up we conclude that
	\begin{align}
		(2\pi)^n\widehat{\mu}_{(2\pi)^{-n},\infty}\geq (2\pi)^n\widehat{\mu}_{c_1,\infty}+(2\pi)^n\widehat{\mu}_{(2\pi)^{-n}-c_1,\infty}
		> (2\pi)^{n}\widehat{\mu}_{(2\pi)^{-n},\infty},\label{eq:4.23}
	\end{align}
	a contradiction. Thus $\|v\|_{L_x^2}^2=(2\pi)^{-n}$. Moreover, if $\widehat{E}^\infty(v)>\widehat{\mu}_{(2\pi)^{-n},\infty}$, then repeating the previous calculation, we arrive at the contradiction \eqref{eq:4.23} again. Thus $v$ is a minimizer of $\widehat{\mu}_{(2\pi)^{-n},\infty}$, and $u_\tau\to v$  in $H_{x,y}^1$  as $\tau\to\infty$. This completes the desired proof.
\end{proof}

\begin{lemma}\label{lem:4.5}
	Let  $\alpha\in(0,\frac{4}{d+n})$. Then there exists some $\tau_*>0$ such that for all $\tau>\tau_*$ we have $\nabla_y u_\tau=0$.
\end{lemma}

\begin{proof}
	From Proposition \ref{lem:4.4}, we can assume that $u_\tau\to v$  in $H_{x,y}^1$.
	Define $w_\tau:=\sqrt{-\Delta_y}u_\tau$. Then applying $\sqrt{-\Delta_y}$ to \eqref{eq:4.16} we obtain
	\begin{align}
		&-\left( a+\frac{b}{\tau^{\frac{4}{\alpha}+2-d}} \|\nabla_{x}u_\tau\|_2^2+\frac{b}{2\tau^{\frac{4}{\alpha}-d}}\|\nabla_{y}u_\tau\|^2_2\right) \Delta_{x}w_\tau\nonumber\\
		&-\left( a\tau^2+\frac{b}{2\tau^{\frac{4}{\alpha}-d}}\|\nabla_{x}u_\tau\|^2_2+b\tau^{2+d-\frac{4}{\alpha}}\|\nabla_{y}u_\tau\|^2_2\right) \Delta_{y}w_\tau +\omega_\tau w_\tau=\sqrt{-\Delta_y}(|u_\tau|^\alpha u_\tau).\label{eq:4.24}
	\end{align}
	Then testing \eqref{eq:4.24} with $w_\tau$ yields
	\begin{align*}
		&a\tau^2\|\nabla_y w_\tau\|_2^2 +\frac{b}{2\tau^{\frac{4}{\alpha}-d}}\|\nabla_{x}u_\tau\|^2_2\|\nabla_{y}w_\tau\|^2_2+ b\tau^{2+d-\frac{4}{\alpha}}\|\nabla_y u_\tau\|_2^2 \|\nabla_y w_\tau\|_2^2  +\frac{b}{2\tau^{\frac{4}{\alpha}-d}}\|\nabla_{y}u_\tau\|^2_2\|\nabla_{x}w_\tau\|^2_2\\
		&\quad +a\|\nabla_x w_\tau\|_2^2	+ \frac{b}{\tau^{\frac{4}{\alpha}+2-d}}\|\nabla_x u_\tau\|_2^2\|\nabla_x w_\tau\|_2^2+\omega_\tau\|w_\tau\|_2^2 -\int_{\mathbb{R}^d\times\mathbb{T}^n}
		\sqrt{-\Delta_y}(|u_\tau|^\alpha u_\tau)w_\tau\,dxdy=0.
	\end{align*}
	In turn this gives
	\begin{align}
		0=&	\left[ a(\tau^2-1)\|\nabla_y w_\tau\|_2^2 +\frac{b}{2\tau^{\frac{4}{\alpha}-d}}\|\nabla_{x}u_\tau\|^2_2\|\nabla_{y}w_\tau\|^2_2+ b\tau^{2+d-\frac{4}{\alpha}}\|\nabla_y u_\tau\|_2^2 \|\nabla_y w_\tau\|_2^2 +\frac{b}{2\tau^{\frac{4}{\alpha}-d}}\|\nabla_{y}u_\tau\|^2_2\|\nabla_{x}w_\tau\|^2_2\right. \nonumber\\ 
		& \left. + \frac{b}{\tau^{\frac{4}{\alpha}+2-d}}\|\nabla_x u_\tau\|_2^2\|\nabla_x w_\tau\|_2^2  -(\alpha+1)\int_{\mathbb{R}^d\times\mathbb{T}^n}
		\sqrt{-\Delta_y}(|v|^\alpha u_\tau)w_\tau\,dxdy\right] \nonumber\\
		&+ 	a\|\nabla_y w_\tau\|_2^2  + a\|\nabla_x w_\tau\|_2^2  + \hat{\omega}\|w_\tau\|_2^2 +\int_{\mathbb{R}^d\times\mathbb{T}^n}
		\sqrt{-\Delta_y}(u_\tau((\alpha+1)|v|^\alpha-|u_\tau|^\alpha))w_\tau\,dxdy\nonumber\\
		& +(\omega_\tau-\hat{\omega})\|w_\tau\|_2^2\nonumber\\
		=&I+II+III.\nonumber
	\end{align}
	Using standard elliptic regularity theory and bootstrap arguments \cite[Theorem 8.1.1]{BookCazenave2003} and Sobolev embedding we have $v\in L^\infty(\mathbb{R}^d)$ ($v$ satisfies \eqref{eq:4.20}). Next we fix an orthonormal basis of eigenfunctions for $-\Delta_y$ (see \cite{Book2003Stein}), that is, $-\Delta_y\varphi_k=\lambda_k\varphi_k$ and $\varphi_0 = const$. Then we can write the following development:
	\begin{equation*}
		w_\tau(x,y)=\sum_{k\in\mathbb{N}\backslash\{0\}}\sigma_{\tau,k}(x)\varphi_k(y).
	\end{equation*}
	Summing up, we conclude that
	\begin{align}\label{eq:4.26}
		I&\geq a(\tau^2-1)\|\nabla_y w_\tau\|_2^2-(\alpha+1)\int_{\mathbb{R}^d\times\mathbb{T}^n}
		\sqrt{-\Delta_y}(|v|^\alpha u_\tau)w_\tau\,dxdy\nonumber\\ 
		&\geq a(\tau^2-1)\sum_{k\neq0}|\lambda_k|^2\int_{\R^d}|\sigma_{\tau,k}(x)|^2\,dx-C(\|v\|_{L_x^\infty},\alpha)\sum_{k\neq0}\int_{\R^d}|\sigma_{\tau,k}(x)|^2\,dx\geq 0
	\end{align}
	for all sufficiently large $\tau$. Moreover, similar to the proof of \cite[Lemma 3.6]{Tzvetkov2014AnalPDE}, by
	using Proposition \ref{lem:4.4},  H\"older inequality and
	the Sobolev embedding $ H^1_{x,y}\hookrightarrow L^{\alpha +2}_{x,y}$,  we infer that
	\begin{align}\label{eq:4.27}
		\left| \int_{\mathbb{R}^d\times\mathbb{T}^n}
		\sqrt{-\Delta_y}(u_\tau((\alpha+1)|v|^\alpha-|u_\tau|^\alpha))w_\tau\,dxdy\right| \leq o_\tau(1)\|w_\tau\|^2_{H_{x,y}^1}.
	\end{align}
	Therefore,  \eqref{eq:4.27} and the fact $\omega_\tau=\hat{\omega}+o_\tau(1)$  imply
	\begin{align}\label{eq:4.28}
		II+III\geq \|w_\tau\|_{H_{x,y}^1}^2(1-o_\tau(1))\geq C \|w_\tau\|_{H_{x,y}^1}^2,
	\end{align}
	for all $\tau>\tau_0$ with some sufficiently large $\tau_0$. By combining   \eqref{eq:4.26} and \eqref{eq:4.28}, we deduce $w_\tau = 0$ for $\tau$ large enough and the desired proof is complete.
\end{proof}
In the following, we prove Theorem \ref{th:1.2} by combining Lemmas \ref{lem:4.1} and \ref{lem:4.5}.
\begin{proof}[Proof of Theorem \ref{th:1.2}]	
	\textbf{Case (i):}
	When  $\alpha\in(0,\frac{4}{d+n})$, Lemma \ref{lem:4.5} yields the existence of $\tau_* \in (0,\infty)$ such that  any minimizer $u_\tau$ of $\mu_{1,\tau}$ satisfies $\nabla_y u_\tau = 0$ for all $\tau \in (\tau_*,\infty)$. 
	For $\kappa>0$ and $c>0$ define the operator $T_1^\kappa$ by
	$$
	T_1^\kappa u(x,y):=\kappa^{\frac{2}{\alpha}}u(\kappa x,y).
	$$
	Let $\kappa_c:=c^{\frac{1}{d-\frac{4}{\alpha}}}$, then $u\mapsto T_1^{\kappa_c} u$ defines a bijection between $\mathcal{S}(c)$ and $\mathcal{S}(1)$. Direct calculation also shows that 
	$$
	E(u)=c^{\frac{d-2-4/\alpha}{d-4/\alpha}}E^{\tau}(T_1^{\kappa_c} u)\quad\text{for}\quad\tau = \kappa_c,
	$$
	thus $m_{c}=c^{\frac{d-2-4/\alpha}{d-4/\alpha}}\mu_{1,\tau}$ for $\tau = \kappa_c$.
	Moreover, by same rescaling arguments we have 
	$\widehat{m}_{(2\pi)^{-n}c}=\widehat{\mu}_{(2\pi)^{-n}c,1}=c^{\frac{d-2-4/\alpha}{d-4/\alpha}}\widehat{\mu}_{(2\pi)^{-n},\tau}$
	for $\tau = \kappa_c$. 
Notice  that the mapping $c\mapsto \kappa_c$ is strictly monotone decreasing on $(0,\infty)$. Thus, if $u_c$ be a minimizer of $m_c$, we know that there exists some $c_{*}\in(0,\infty)$  such that for any
 $c\in(0,c_{*})$ we have $\nabla_y u_c= 0$. 
	
	\textbf{Case (ii):}	 When $\alpha\in(2,\min\{\frac{8}{d+n},\frac{4}{d}\})$, Lemma \ref{lem:4.1} yields the existence of $\lambda_* \in (0,\infty)$ such that for all $\lambda \in (0,\lambda_*)$ any minimizer $u_\lambda$ of $m_{1,\lambda}$ satisfies $m_{1,\lambda}<(2\pi)^n \widehat{m}_{(2\pi)^{-n},\lambda}$. This implies that for all $\lambda \in (0,\lambda_*)$ we have $\nabla_y u_\lambda \neq 0$. In fact, if $\nabla_yu_\lambda=0$, then 
	$$m_{1,\lambda}=(2\pi)^n\left( \frac{a\lambda^{\frac{4-d}{\alpha-2}+1-\frac{d}{2}}}{2}\|\nabla_x u\|_{L_x^2}^2+\frac{b(2\pi)^n}{4}\|\nabla_x u\|_{L_x^2}^4
	-\frac{1}{\alpha+2}\|u\|_{L_x^{\alpha+2}}^{\alpha+2}\right) >(2\pi)^n\widehat{m}_{(2\pi)^{-n},\lambda}.$$
	Moreover, 	for $\kappa>0$	and $c>0$,   define the operator $T_2^\kappa$ by
	$$
	T_2^\kappa u(x,y):=\kappa^{\frac{4-d}{\alpha-2}}u(\kappa x,y).
	$$
	Let $\kappa_c=c^{\frac{1}{d-\frac{8-2d}{\alpha-2}}}$, then $u\mapsto T_2^{\kappa_c} u$ defines a bijection between $\mathcal{S}(c)$ and $\mathcal{S}(1)$. Direct calculation  shows that
	$$
	E(u)=c^{\frac{(4-d)(\alpha+2)/(\alpha-2)-d}{d-(8-2d)/(\alpha-2)}}E_\lambda(T_2^{\kappa_c} u)\quad \text{ for } \quad\lambda= \kappa_c^{2},
	$$
	thus 
	$$
	m_{c}=c^{\frac{d-(4-d)(\alpha+2)/(\alpha-2)}{d-(8-2d)/(\alpha-2)}}m_{1,\lambda}
	$$ 
	for $\lambda=\kappa_c^{2}$. By same rescaling arguments we also infer that
	\begin{align}\label{invariant hat minimizing}
	\widehat{m}_{(2\pi)^{-n}c}=	\widehat{m}_{(2\pi)^{-n}c,1}=c^{d-\frac{(4-d)(\alpha+2)/(\alpha-2)}{d-(8-2d)/(\alpha-2)}}\widehat{m}_{(2\pi)^{-n},\lambda}\quad\text{for $\lambda =\kappa_c^{2}$.}
	\end{align}
Notice also that the mapping $c\mapsto\kappa_c$ is strictly monotone decreasing on $(0,\infty)$. Thus we know that there exists some $c_{**}\in(0,\infty)$  such that for any $c\in(c_{**},\infty)$ we have $m_{c}<(2\pi)^n\widehat{m}_{(2\pi)^{-n}c,1}=(2\pi)^n\widehat{m}_{(2\pi)^{-n}c}$, and then for $u_c$ be a minimizer of $m_c$ we infer that $\nabla_y u_c\neq 0$  with  $c\in(c_{**},\infty)$. 
\end{proof}

\section{Existence and $y$-dependence of normalized ground states: Mass supercritical case}

\subsection{Existence  of normalized ground states}\label{sub2.1}
In this subsection, we study the existence of normalized solutions of problem \eqref{eq:1.1} in the mass supercritical. First, We shall also frequently use results arising from the variational problem $\widehat{M}_c$, which we summarize in the following. The results are classical and we refer e.g. to \cite{BookCazenave2003,Jeanjean1997NA,Ye2015MMAS} for details of the corresponding proofs.

\begin{proposition}\label{lem:1.5}
	The following statements hold true:
	\begin{itemize}
		\item[(i)]For any $c>0$,  $Q_c\in \widehat{\mathcal{S}}(c)$ is an minimizer of $\widehat{M}_c$ (defined in \eqref{hMc}) and satisfies the following equation
		\begin{align}\label{eqrd}
			-\left( a+b\int_{\mathbb{R}^d}|Q_c|^2dx\right) \Delta_x Q_c+\omega_c Q_c=|Q_c|^\alpha Q_c
		\end{align}
		with some $\omega_c>0$.
		\item[(ii)] Any solution $Q_c\in H^1(\mathbb{R}^d)$ of \eqref{eqrd} with $\omega_c>0$ is of class $W^{3,p}(\mathbb{R}^d)$ for all $p\in[2,\infty)$.
		\item[(iii)] Any solution $Q_c\in H^1(\mathbb{R}^d)$ of \eqref{eqrd} satisfies $\widehat{P}(Q_c)=0$.
		\item[(iv)] The mapping $c\mapsto \widehat{M}_c$ is strictly monotone decreasing and continuous on $(0,\infty)$.
	\end{itemize}
\end{proposition}
Now, we consider some useful conclusions about $M_c$  and $P(u)$.
\begin{lemma}\label{lem:5.1}
	Assume that $\alpha\in(8,\infty)$, then  we have $M_c\in(0,\infty)$ for any $c\in(0,\infty)$, here $M_c$ is defined in \eqref{eq:Mc}.
\end{lemma}
\begin{proof}
	That $M_c<\infty$ follows from $V(c)\neq\varnothing$. Next, let  $\{u_k\}\subset V(c)$ such that $E(u_k)=M_c+o_k(1)$. Then
	\begin{align*}
		E(u_k)&=E(u_k)-\frac{2}{\alpha }P(u_k)\\
		&=\frac{a}{2}\|\pt_y u_k\|_2^2+\frac{b}{4}\|\pt_y u_k\|_2^4+\left( \frac{1}{2}-\frac{2}{\alpha }\right)\left( a+b\|\pt_y u_k\|_2^2\right)  \|\pt_x u_k\|_2^2+b\left( \frac{1}{4}-\frac{2}{\alpha }\right) \|\pt_x u_k\|_2^4.
	\end{align*}
	Combining with $\alpha>8$  we infer that $\{u_k\}$ is a bounded sequence in $H_{x,y}^1$. Now by the Gagliardo-Nirenberg inequality \eqref{eq:gn2}  we deduce
	\begin{align*}
		a \|\pt_x u_k\|_2^2+b \|\pt_x u_k\|_2^4&\leq\frac{\alpha}{2(\alpha+2)}\|u_k\|_{\alpha+2}^{\alpha+2}\\
		&\leq C(\alpha) \|\pt_x u_k\|_2^{\frac{\alpha }{2}}\|u_k\|_2^{2}
		\left( \| u_k\|_{2}^{\frac{\alpha}{2}}+\|\pt_y u_k\|_{2}^{\frac{\alpha}{2}}\right) \leq C \|\pt_x u_k\|_2^{\frac{\alpha }{2}},
	\end{align*}
	this implies 
	\begin{align}
		\liminf_{k\to\infty}\|u_k\|_{\alpha+2}^{\alpha+2}\geq C\liminf_{k\to\infty}\left( a \|\pt_x u_k\|_2^2+b \|\pt_x u_k\|_2^4\right) >0.\label{eq:5.1}
	\end{align}
	Thus, we obtain
	\begin{align*}
		M_c=\lim_{k\to\infty}E(u_k)&=\lim_{k\to\infty}E(u_k)-\frac{2}{\alpha }P(u_k)\geq \left( \frac{1}{4}-\frac{2}{\alpha}\right) \liminf_{k\to\infty}\left( a \|\pt_x u_k\|_2^2+b \|\pt_x u_k\|_2^4\right)>0,
	\end{align*}
	which completes the proof.
\end{proof}
In the following, we prove some useful properties of the mappings $t\mapsto P(u\ast t)$ and $c\mapsto M_c$, where $u\ast t$  defined in \eqref{eq:L2scaling}.
\begin{lemma}\label{lem:5.2}
	Let $c>0$ and $u\in \mathcal{S}(c)$. Then the following statements hold true:
	\begin{enumerate}
		\item[(1)] $\frac{\partial}{\partial t}E(u\ast t)=t^{-1} P(u\ast t)$ for all $t>0$;
		\item[(2)] there exists some $t^*=t^*(u)>0$ such that $u\ast {t^*}\in V(c)$;
		\item[(3)] $P(u\ast t)>0$ when $t\in(0,t^*)$ and  $P(u\ast t)<0$ when $t\in(t^*,\infty)$;
		\item[(4)] we have $t^*<1$ if and only if $P(u)<0$. Moreover, $t^*=1$ if and only if $P(u)=0$;
		\item[(5)] $E(u\ast t)<E(u\ast {t^*})$ for all $t>0$ with $t\neq t^*$.
	\end{enumerate}
\end{lemma}
\begin{proof}
	Define $f(t):= \frac{\partial}{\partial t}E(u\ast t)$. Then
	\begin{align*}
		f(t)=t(a+b\|\nabla_y u\|_2^2)\|\pt_x u\|_2^2+bt^3\|\pt_x u\|_2^4-\frac{\alpha}{2(\alpha+2)}t^{\frac{\alpha }{2}-1}\|u\|_{\alpha+2}^{\alpha+2},
	\end{align*}
	and
	\begin{align*}
		f'(t)&=(a+b\|\nabla_y u\|_2^2)\|\pt_x u\|_2^2+bt^2\|\pt_x u\|_2^4-\frac{\alpha (\alpha -2)}{4(\alpha+2)}t^{\frac{\alpha }{2}-2}\|u\|_{\alpha+2}^{\alpha+2}.
	\end{align*}
	Using $\alpha>8$ we infer that there exists $t_*$ such that  $f'(t)$ is strictly strictly increasing on $(0,t_*)$ and strictly  decreasing on $(t_*,\infty)$. Since $f'(0)=(a+b\|\nabla_y u\|_2^2)\|\pt_x u\|_2^2>0$ and $f'(t)\to -\infty$ as $t\to\infty$, there exists a $t_0>0$ such that $f'(t)$ is positive on $(0,t_0)$ and negative on $(t_0,\infty)$. Consequently, we conclude that $f(t)$ has a zero at $t^*>t_0$ such that $f(t)$ is positive on $(0,t^*)$ and negative on $(t^*,\infty)$. 
	
	Then, (1) follows  from direct calculation.
	(2) and (3)  follow from the fact (1).
	For (4), let $P(u)<0$. Then
	$$0>P(u)=\frac{P(u^1)}{1}=f(1),$$
	which is only possible as long as $t^*<1$. Conversely, assume that $t^*<1$, then using the fact that $f(t)$ is monotone decreasing on $(t^*,\infty)$ we obtain
	$$P(u)=f(1)<f(t^*)<0. $$
	This completes the proof of (4). To see (5), integration by parts yields
	$$ E(u\ast{t^*})=E(u\ast t)+\int_t^{t^*}f(s)\,ds.$$
	Then (5) follows from the fact that $f(t)$ is positive on $(0,t^*)$ and $f(t)$ is negative on $(t^*,\infty)$.
\end{proof}

\begin{lemma}\label{lem:5.3}
	The mapping $c\mapsto M_c$ is continuous and non-increasing on $(0,\infty)$.
\end{lemma}

\begin{proof}
	Defined 
	\begin{align*}
		g(\beta_1,\beta_2,\beta_3):=\max_{s>0}\{\beta_1s^2+\beta_2s^4-\beta_3s^{\frac{\alpha}{2}}\}
	\end{align*}
	where $\alpha>8$ and $ \beta_1 > 0$, $\beta_2>0$, $\beta_3> 0$ which are totally independent of $t$. In view of \cite{Jean-Luo2013PLMS} the function $g$ is continuous on $\mathbb{R}^+\times\mathbb{R}^+\times\mathbb{R}^+$. 
	
	We now show the monotonicity of $c\mapsto M_c$. It suffices to show that for any $0<c_1<c_2<\infty$ and $\vare>0$ we have
	\begin{align*}
		M_{c_2}\leq M_{c_1}+\vare.
	\end{align*}
	By the definition of $M_{c_1}$ there exists some $u_1\in V(c_1)$ such that
	\begin{align}\label{eq:5.2}
		E(u_1)\leq M_{c_1}+\frac{\vare}{2}.
	\end{align}
	Let $\eta\in C^{\infty}_0(\mathbb{R};[0,1])$ be a cut-off function such that $\eta=1$ for $|x|\leq 1$ and $\eta=0$ for $|x|\geq 2$. For $\delta>0$, define
	\begin{equation*}
		u_{1,\delta}(x,y):= \eta(\delta x)\cdot u_1(x,y).
	\end{equation*}
	Using dominated convergence theorem it is easy to verify that $u_{1,\delta}\to u_1$ in $H_{x,y}^1$ as $\delta\to 0$. Therefore,
	\begin{align*}
		\|\nabla_{x,y}u_{1,\delta}\|_2\to \|\nabla_{x,y} u_1\|_2\quad\text{and}\quad
		\|u_{1,\delta}\|_{\alpha+2}\to \| u_1\|_{\alpha+2},
	\end{align*}
	for all $\alpha>8$ as $\delta\to 0$. Combining with the continuity of $g$ we conclude that
	\begin{align}\label{eq:5.3}
		\max_{t>0}E(u_{1,\delta}\ast t)
		&=\max_{t>0}\left\lbrace \frac{t^2}{2}(a+b\|\pt_{y}u_{1,\delta}\|_2^2)\|\pt_xu_{1,\delta}\|_2^2+\frac{bt^4}{4}\|\pt_xu_{1,\delta}\|_2^4-\frac{t^{\frac{\alpha }{2}}}{\alpha+2}\|u_{1,\delta}\|_{\alpha+2}^{\alpha+2}
		\right\rbrace \nonumber\\
		&\quad+\frac{a}{2}\|\pt_y u_{1,\delta}\|_2^2+\frac{b}{4}\|\pt_y u_{1,\delta}\|_2^4\nonumber\\
		&\leq\max_{s>0}\left\lbrace \frac{t^2}{2}(a+b\|\pt_{y}u_{1}\|_2^2)\|\pt_xu_{1}\|_2^2+\frac{bt^4}{4}\|\pt_xu_{1}\|_2^4-\frac{t^{\frac{\alpha  }{2}}}{\alpha+2}\|u_{1}\|_{\alpha+2}^{\alpha+2}
		\right\rbrace \nonumber\\
		&\quad+\frac{a}{2}\|\pt_y u_{1}\|_2^2+\frac{b}{4}\|\pt_y u_{1}\|_2^4+\frac{\vare}{4}\nonumber\\
		&=\max_{t>0}E( u_1\ast t)+\frac{\vare}{4},
	\end{align}
	for sufficiently small $\delta>0$. Now let $v\in C_0^\infty(\mathbb{R})$ with $\mathrm{supp}\,v\subset
	B(0,4\delta^{-1}+1)\backslash B(0,4\delta^{-1})$ and define
	\begin{equation*}
		v_0:= \frac{(c_2-\|u_{1,\delta}\|_2^2)^{\frac{1}{2}}}{\|v\|_2}\,v.
	\end{equation*}
	Notice that $v_0$ and $u_{1,\delta}$ have compact supports, which also implies $\|v_0\|_2^2=c_2-\|u_{1,\delta}\|_2^2$. Define
	\begin{align*}
		w_\ld:=u_{1,\delta}+v_0\ast \ld,
	\end{align*}
	with some to be determined $\ld>0$. Then
	\begin{align*}
		\|w_\ld\|^p_p=\|\tilde{u}_{1,\delta}\|^p_p+\| v_0\ast\ld\|^p_p,
	\end{align*}
	for all $p\in [2,\infty)$. Particularly, $\|w_\lambda\|_2^2=c_2$. Since $v_0$ is independent of $y\in\T$, we also infer that
	\begin{align*}
		\|\pt_x w_\ld\|_2\to\|\pt_x u_{1,\delta}\|_2,\quad	\|\pt_y w_\ld\|_2=\|\pt_y u_{1,\delta}\|_2\quad\text{and}\quad
		\| w_\ld\|_p\to\| u_{1,\delta}\|_p,
	\end{align*}
	for all $p\in(2,\infty)$ as $\ld\to 0$. Using the continuity of $g$ once again we obtain
	\begin{align*}
		\max_{t>0}E(w_\ld\ast t)\leq \max_{t>0}E(u_{1,\delta}\ast t)+\frac{\vare}{4},
	\end{align*}
	for sufficiently small $\ld>0$. Finally, combing with \eqref{eq:5.2} and \eqref{eq:5.3} we infer that
	\begin{align*}
		M_{c_2}\leq \max_{t>0}E(w_\lambda\ast t)\leq \max_{t>0}E(u_{1,\delta}\ast t)+\frac{\varepsilon}{4}\leq
		\max_{t>0}E(u_1\ast t)+\frac{\varepsilon}{2}=E(u_1)+\frac{\varepsilon}{2}\leq M_{c_1}+\varepsilon,
	\end{align*}
	which implies  $c\mapsto M_c$ is non-increasing on $(0,\infty)$.
	
	Next, we show the continuity of  $c\mapsto M_c$. Since $c\mapsto M_c$ is non-increasing, it suffices to show that for any $c\in(0,\infty)$  we have
	\begin{align}\label{eq:5.4}
		M_c\leq \lim_{n\to \infty}M_{c_k}\quad\text{any sequence }\quad c_k\searrow c.
	\end{align}
	
	Let $\vare>0$ be an arbitrary positive number. For any sequence $c_k\searrow c$, by the definition of $M_{c_k}$ we can find some $u_k\in V(c_k)$ such that
	\begin{align*}
		E(u_k)\leq M_{c_k}+\frac{\vare}{2}\leq M_c+\frac{\vare}{2}.
	\end{align*}
	We define $\tilde{u}_k=(c_k^{-1}c)^{\frac{1}{2}} \cdot u_k:=\rho_k u_k$, then $\|\tilde{u}_k\|_2^2=c$ and $\rho_k\nearrow 1$. Since $u_k\in V(c_k)$  we obtain
	\begin{align*}
		M_{c}+\frac{\vare}{2}&\geq M_{c_k}+\frac{\vare}{2}\geq E(u_k)-\frac{2}{\alpha}P(u_k)\\
		&=\frac{a}{2}\|\pt_y u_k\|_2^2+\frac{b}{4}\|\pt_y u_k\|_2^4+\left( \frac{1}{2}-\frac{2}{\alpha}\right)(a+b\|\pt_{y}u_{n}\|_2^2) \|\pt_x u_k\|_2^2+b\left( \frac{1}{4}-\frac{2}{\alpha}\right) \|\pt_x u_k\|_2^4.
	\end{align*}
	Thus $\{u_k\}$ is bounded in $H^1_{x,y}$ and up to a subsequence we infer that there exist $A,B,C\geq 0$ such that
	\begin{align*}
		\|\pt_xu_k\|_2^2=A+o_k(1),\quad\|\pt_y u_k\|_2^2=B+o_k(1),\quad\|u_k\|_{\alpha+2}^{\alpha+2}=C+o_k(1).
	\end{align*}
	Arguing as in the proof of Lemma \ref{lem:5.1}, we may use the fact $P(u_k)=0$ and Lemma \ref{lem:gn2} to deduce $A,C>0$ and by previous arguments we know that $g$ is continuous at the point $(\frac{(a+bB)A}{2},\frac{bA^2}{4},\frac{C}{\alpha+2})$. Using also the fact that $\rho_k\nearrow 1$ we conclude that
	\begin{align*}
		M_{c}\leq \max_{t>0}E(\tilde{u}_k\ast t)
		=&\max_{t>0}\left\lbrace \frac{at^2\rho_k^2}{2}\|\pt_x{u}_{k}\|_2^2+\frac{bt^4\rho_k^4}{4}\|\pt_x{u}_{k}\|_2^2(\|\pt_x{u}_{k}\|_2^2+\|\pt_{y}u_{k}\|_2^2)
		-\frac{t^{\frac{\alpha}{2}}\rho_k^{\alpha+2}}{\alpha+2}\|{u}_{k}\|_{\alpha+2}^{\alpha+2}\right\rbrace \nonumber\\
		&+\frac{a\rho_k^2}{2}\|\pt_y u_k\|_2^2+\frac{b\rho_k^4}{4}\|\pt_y u_k\|_2^4\nonumber\\
		\leq&\max_{t>0}\left\lbrace \frac{t^2(a+bB)A}{2}+\frac{t^4bA^2}{4}
		- \frac{t^{\frac{\alpha}{2}}C}{\alpha+2}\right\rbrace 
		+\frac{a}{2}\|\pt_y u_k\|_2^2+\frac{b}{4}\|\pt_y u_k\|_2^4+\frac{\vare}{4}\nonumber\\
		\leq& \max_{t>0}\left\lbrace \frac{at^2}{2}\|\pt_x u_k\|_2^2+\frac{bt^2}{4}\|\pt_x u_k\|_2^2(\|\pt_x{u}_{n}\|_2^2+\|\pt_{y}u_{n}\|_2^2)
		- \frac{t^{\frac{\alpha}{2}}}{\alpha+2}\| u_k\|_{\alpha+2}^{\alpha+2}\right\rbrace \\
		&+\frac{a}{2}\|\pt_y u_k\|_2^2+\frac{b}{4}\|\pt_y u_k\|_2^4+\frac{\vare}{2}\nonumber\\
		=&\max_{t>0}E(u_k\ast t)+\frac{\vare}{2}=E(u_k)+\frac{\vare}{2}\leq M_{c_k}+\vare,
	\end{align*}
	by choosing $n$ sufficiently large.  The continuity follows from the arbitrariness of $\vare$. 
\end{proof}
Since in general the prefactors of $-\pt_x u_c$ and $-\pt_y^2 u_c$ do not coincide in the Lagrange multiplier equation of $u_c$, we do not know whether a simple proof via the Pohozaev identity is applicable. In this case, we appeal to a subtle deformation argument (see for instance \cite{Jean-Luo2013PLMS,willem1997}) to overcome this difficulty.
We shall firstly introduce some preliminary concepts given by \cite{Lions1983ARMAii}. First recall that $\mathcal{S}(c)$ is a submanifold of $H_{x,y}^1$ with codimension $1$. Moreover, the tangent space $T_u\mathcal{S}(c)$ for a point $u\in \mathcal{S}(c)$ is given by
$$ T_u \mathcal{S}(c)=\left\lbrace v\in H_{x,y}^1:\la u,v\ra_{L_{x,y}^2}=0\right\rbrace .$$
We also denote the tangent bundle of $\mathcal{S}(c)$ by $T\mathcal{S}(c)$. Next, the energy functional $E|_{\mathcal{S}(c)}$  is a $C^1$-functional and for any $u\in \mathcal{S}(c)$ and $v\in T_u\mathcal{S}(c)$ we have
\begin{align*}
	\la E'|_{\mathcal{S}(c)}(u),v\ra =\la E'(u),v\ra .
\end{align*}
We use $\|E'|_{\mathcal{S}(c)}(u)\|_*$ to denote the dual norm of $E'|_{\mathcal{S}(c)}(u)$ in the cotangent space $(T_u\mathcal{S}(c))^*$, i.e.
$$\|E'|_{\mathcal{S}(c)}(u)\|_*:=\sup_{v\in T_u \mathcal{S}(c),\,\|v\|_{H^1_{x,y}}\leq 1}|\la E'|_{\mathcal{S}(c)}(u),v\ra|. $$
\begin{definition}\label{definiton of mp geometry}
	We say that $E(u)$ has a mountain pass geometry on $\mathcal{S}(c)$ at the level $\gamma_c$ if there exists some $\rho>0$ and $\vare\in\left( 0,\min\left\lbrace M_c,\frac{\alpha-8}{3\alpha}\right\rbrace \right) $ such that
	\begin{align*}
		\gamma_c:=\inf_{h\in\Gamma(c)}\max_{t\in[0,1]}E(h(t))>\max\left\lbrace \sup_{h\in\Gamma(c)}E(h(0)),\sup_{h\in\Gamma(c)}E(h(1))\right\rbrace ,
	\end{align*}
	where
	\begin{align*}
		\Gamma(c):=\{h\in C([0,1];\mathcal{S}(c)):h(0)\in A_{\rho,\varepsilon},E(h(1))<0\}
	\end{align*}
	and
	\begin{align*}
		A_{\rho,\vare}:=\{u\in \mathcal{S}(c):\|\pt_x u\|_2^2<\rho,\quad 2a\|\pt_y u\|_2^2+b\|\pt_y u\|_2^4 \leq 4(M_c-\vare)\}.
	\end{align*}
\end{definition}

We next show that one obtains mountain geometry by choosing $k$ and $\vare$ sufficiently small.

\begin{lemma}\label{lem:5.5}
	There exist $\rho>0$ and $\vare\in(0,M_c)$ such that
$E(u)$ has a mountain pass geometry on $\mathcal{S}(c)$ at the level $\gamma_c$ in the sense of Definition \ref{definiton of mp geometry}, and $\gamma_c=M_c$. 
\end{lemma}

\begin{proof}
	We firstly prove that by choosing $\rho$ sufficiently small we have $P(u)>0$ for all $u\in A_{\rho,\vare}$, where $\rho$ is independent of the choice of  $\vare\in (0,M_c)$. Indeed, by Lemma \ref{lem:gn2} we obtain
	\begin{align*}
		P(u)=(a+b\|\pt_{y}u\|_2^2) \|\pt_x u\|_2^2+b \|\pt_x u\|_2^4-\frac{\alpha}{2(\alpha+2)}\|u\|_{\alpha+2}^{\alpha+2}
		\geq b \|\pt_x u\|_2^4-C\|\pt_x u\|_2^{\frac{\alpha}{2}}>0,
	\end{align*}
	for some sufficiently small $\rho$. Next we choose the number $\vare$. Arguing as in \eqref{eq:5.1} we know that there exists some $\beta=\beta(c)>0$ such that if $\{u_k\}\subset V(c)$ is a minimizing sequence for $M_c$, then
	\begin{align}
		\liminf_{k\to\infty}\left[ a\left( \frac{1}{2}-\frac{2}{\alpha}\right) \|\pt_x u_k\|_2^2+b\left( \frac{1}{4}-\frac{2}{\alpha}\right) \|\pt_x u_k\|_2^4\right] \geq \beta.\label{eq:5.6}
	\end{align}
	We may shrink $\beta$ further such that $\beta<4M_c$. Hence for any minimizing sequence $\{u_k\}$ we  have
	\begin{align}\label{eq:5.7}
		M_c+\frac{\beta}{4}&\geq \frac{a}{2}\|\pt_y u_k\|_2^2+\frac{b}{4}\|\pt_y u_k\|_2^4+\left( \frac{1}{2}-\frac{2}{\alpha}\right)\left( a \|\pt_x u_k\|_2^2+b \|\pt_x u_k\|_2^2\|\pt_{y}u_{k}\|_2^2\right) +b\left( \frac{1}{4}-\frac{2}{\alpha}\right) \|\pt_x u_k\|_2^4\nonumber\\
		&>\frac{a}{2}\|\pt_y u_k\|_2^2+\frac{b}{4}\|\pt_y u_k\|_2^4+\frac{\beta}{2},
	\end{align}
	for all sufficiently large $k$. Now we set $\vare=\frac{\beta}{4}$. Therefore for $u\in A_{\rho,\vare}$, using also Lemma \ref{lem:gn2} we infer that
	\begin{align*}
		E(u)&\leq \frac{a}{2}\|\pt_y u\|_2^2+\frac{b}{4}\|\pt_y u\|_2^2+\left( \frac{a}{2}+\frac{b}{2}\|\pt_y u\|_2^2\right)  \|\pt_x u\|_2^2+\frac{b}{4} \|\pt_x u\|_2^4+C\|\pt_x u\|_2^{\frac{\alpha}{2}}\\
		&	\leq M_c-\vare+\left( \frac{a}{2}+\frac{b}{2}\|\pt_y u\|_2^2\right)\rho+\frac{b}{4}\rho^2+C\rho^{\frac{\alpha}{4}}.
	\end{align*}
	We now choose $\rho=\rho(\vare)$ sufficiently small  (without changing the fact that $P(u)>0$ for all $u\in A_{\rho,\vare}$) such that 
	$$
	\left( \frac{a}{2}+\frac{b}{2}\|\pt_y u\|_2^2\right)\rho+\frac{b}{4}\rho^2+C\rho^{\frac{\alpha}{4}}<\frac{\vare}{2}.
	$$ 
	Thus for this choice of $\rho$ and $\vare$, we have $E(u)<M_c-\frac{\varepsilon}{2}<M_c$ for all $u\in A_{\rho,\vare}$.  	It is easy to check that $\|\pt_y(u\ast t)\|_2^2=\|\pt_y u\|_2^2$ for all $t\in(0,\infty)$ and $\|\pt_x (u\ast t)\|_2^2=t^2\|\pt_x u\|_2^2\to 0$ as $t\to 0$. We then fix some $t_0>0$ sufficiently small such that
	\begin{align*}
		\|\pt_x(u\ast{t_0})\|_2^2\leq\|\pt_xu\|_2^2<\rho\quad\text{and}\quad\frac{a}{2}\|\pt_y (u\ast{t_0})\|_2^2+\frac{b}{4}\|\pt_y (u\ast{t_0})\|_2^4\leq M_c-\vare,
	\end{align*} 
	namely, $u\ast{t_0}\in A_{\rho,\vare}$. In addition,
	\begin{align*}
		E(u\ast t)=\frac{a}{2}\|\pt_y u\|_2^2+\frac{b}{4}\|\pt_y u\|_2^4 +\frac{t^2}{2}(a+b\|\pt_y u\|_2^2)\|\pt_xu\|_2^2+\frac{bt^4}{4}\|\pt_xu\|_2^4-\frac{t^{\frac{\alpha}{2}}}{\alpha+2}\|u\|_{\alpha+2}^{\alpha+2}\to -\infty,
	\end{align*}
	as $t\to\infty$. We then fix some $t_1$ sufficiently large such that $E(u*t_1)<0$. Now define $g\in C([0,1];\mathcal{S}(c))$ by
	\begin{align}\label{eq:5.8}
		g(t):=u\ast({t_0+(t_1-t_0)t})\in \Gamma(c).
	\end{align}
	This is $\Gamma(c)\neq\emptyset$. On the other hand, using the fact that $P(u)>0$ for all $u\in A_{\rho,\vare}$,  we have $P(h(0))>0$ for any $h\in \Gamma(c)$.  Assume the contrary that there exists some $h\in\Gamma(c)$ such that $P(h(1))\geq 0$. Then
	\begin{align*}
		0>E(h(1))&\geq E(h(1))-\frac{2}{\alpha}P(h(1))\\
		&=\frac{a}{2}\|\pt_y h(1)\|_2^2+\frac{b}{4}\|\pt_y h(1)\|_2^4+\left( \frac{1}{2}-\frac{2}{\alpha }\right)\left( a+b\|\pt_y h(1)\|_2^2\right)  \|\pt_x h(1)\|_2^2+b\left( \frac{1}{4}-\frac{2}{\alpha }\right) \|\pt_x h(1)\|_2^4\\
		&\geq a\left( \frac{1}{2}-\frac{2}{\alpha}\right) \|\pt_x h(1)\|_2^2+b\left( \frac{1}{4}-\frac{2}{\alpha}\right) \|\pt_x h(1)\|_2^4\geq 0,
	\end{align*}
	a contradiction. Next, by continuity of $h$ there exists some $t^*\in(0,1)$ such that $P(h(t^*))=0$, which implies  $h(t^*)\in V(c)$. 
	It is easy to $E(u)$ has a mountain pass geometry on $\mathcal{S}(c)$ at the level $\gamma_c$.
	
	Next we show that $\gamma_c=M_c$.
 Let $\{u_k\}\subset V(c)$ be the given minimizing sequence satisfying \eqref{eq:5.6}.  For any $\kappa\in(0,\beta/4)$, set $u:=u_k$ for some  sufficiently large $k\in\N$ such that $E(u)\leq M_c+\kappa$ and 
	$$
	a\left( \frac{1}{2}-\frac{2}{\alpha}\right) \|\pt_x u\|_2^2+b\left( \frac{1}{4}-\frac{2}{\alpha}\right) \|\pt_x u\|_2^4\geq \beta/2.
	$$
	Then by \eqref{eq:5.7} we know that 
	\begin{align*}
		\frac{a}{2}\|\pt_y u\|_2^2+\frac{b}{4}\|\pt_y u\|_2^4\leq M_c-\vare\quad\text{for all} \quad\kappa\in(0,\beta/4).
	\end{align*} 
	 By definition of $\gamma_c$ and $g\in C([0,1];\mathcal{S}(c))$ defined in \eqref{eq:5.8} we have
	\begin{align*}
		\gamma_c\leq \max_{t\in[0,1]}E(g(t))=E(u)\leq M_c+\kappa.
	\end{align*}
	Since $\kappa$ can be chosen arbitrarily small, we conclude that $\gamma_c\leq M_c$. From $h(t)\cap V(c)\neq\emptyset$ for all $h(t)\in\Gamma(c)$, we have
	\begin{align*}
		\max_{t\in[0,1]}E(h(t))\geq M_c,
	\end{align*}
	then we deduce $\gamma_c\geq M_c$, which completes the  proof.
\end{proof}

Let now
$$\tilde{S}(c):=\{u\in \mathcal{S}(c):E'|_{\mathcal{S}(c)}(u)\neq 0\}.$$
According to \cite[Lemma 4]{Lions1983ARMAii} there exists a locally Lipschitz pseudo gradient vector field $Y:\tilde{S}(c)\to T\mathcal{S}(c)$ such that $Y(u)\in T_u\mathcal{S}(c)$,
\begin{align}\label{eq:5.9}
	\|Y(u)\|_{H_{x,y}^1}\leq 2\|E'|_{\mathcal{S}(c)}(u)\|_*\quad\text{and}\quad \la E'|_{\mathcal{S}(c)}(u),Y(u)\ra\geq \|E'|_{\mathcal{S}(c)}(u)\|_*^2
\end{align}
for $u\in \tilde{S}(c)$.
Borrow an idea from the proof of \cite[Lemma 6.1]{Jean-Luo2013PLMS}, we can obtain the following lemma.
\begin{lemma}\label{lem:5.6}
	For any $c>0$, a minimizer $u_c$ of $M_c$ is a solution of problem \eqref{eq:1.1} for some $\omega\in\R$.
\end{lemma}

\begin{proof}
	By Lagrange multiplier  theory we know that $u_c$ solves problem \eqref{eq:1.1} is equivalent to $E'|_{\mathcal{S}(c)}(u_c)=0$. We  assume that $\|E'|_{\mathcal{S}(c)}(u_c)\|_*\neq 0$, which implies that there exists some $\delta>0$ and $\sigma>0$ such that
	\begin{align}\label{eq:5.10}
		v\in B_{u_c}(3\delta)\quad\text{and}\quad\|E'|_{\mathcal{S}(c)}(v)\|_*\geq \sigma,
	\end{align}
	where $B_{u}(\delta):=\{v\in \mathcal{S}(c):\|u-v\|_{H_{x,y}^1}\leq \delta\}$. Let $\rho$ and $\vare$ be given according to Lemma \ref{lem:5.5} ( if necessary such that $\varepsilon<\frac{4\alpha-32}{12\alpha}$). Define
	\begin{align*}
		\vare_1:=\frac{1}{4}\left( M_c-\max\left\lbrace \sup_{h\in\Gamma(c)}E(h(0)),\sup_{h\in\Gamma(c)}E(h(1))\right\rbrace \right) ,
	\end{align*}
	and
	\begin{align*}
		\vare_2:=\min\{\vare_1,M_c/4,\sigma\delta/4\}.
	\end{align*}
	Let the sets $\mathcal{A},\mathcal{B}$  be given by
	\begin{align*}
		\mathcal{A}:= \mathcal{S}(c)\cap E^{-1}([M_c-2\vare_2,M_c+2\vare_2]),
	\end{align*}
	and
	\begin{align*}
		\mathcal{B}:=B_{u_c}(2\delta)\cap E^{-1}([M_c-\vare_2,M_c+\vare_2]).
	\end{align*}
	Moreover, set $f:\mathcal{S}(c)\to [0,\delta]$ satisfies
	\begin{align*}
		f(u):=\frac{\delta\,\mathrm{dist}(u,\mathcal{S}(c)\setminus \mathcal{A})}{\mathrm{dist}(u,\mathcal{S}(c)\setminus \mathcal{A})+\mathrm{dist}(u,\mathcal{B})}=
		\left\{
		\begin{array}{ll}
			\delta,&\text{if $u\in\mathcal{B}$},\\
			0,&\text{if $u\notin \mathcal{A}$}.
		\end{array}
		\right.
	\end{align*}
	Next, we define the pseudo gradient flow $W:\mathcal{S}(c)\to H_{x,y}^1$ by
	\begin{align*}
		W(u):=
		\left\{
		\begin{array}{ll}
			-f(u)\|Y(u)\|^{-1}_{H_{x,y}^1}Y(u),&\text{if $u\in\tilde{S}(c)$},\\
			0,&\text{if $u\in \mathcal{S}(c)\setminus \tilde{S}(c)$}.
		\end{array}
		\right.
	\end{align*}
	One easily verifies that $W$ is a locally Lipschitz function from $\mathcal{S}(c)$ to $H_{x,y}^1$. Then by standard arguments (see \cite[Lem. 6]{Lions1983ARMAii}) there exists a mapping $\eta:\R\times \mathcal{S}(c)\to \mathcal{S}(c)$ such that $\eta(1,\cdot)\in C(\mathcal{S}(c);\mathcal{S}(c))$ and $\eta$ solves the differential equation
	\begin{align*}
		\frac{d}{dt}\eta(t,u)=W(\eta(t,u)),\quad\eta(0,u)=u\quad\text{for any} \quad u\in \mathcal{S}(c).
	\end{align*}
	We now claim that $\eta$ satisfies the following properties:
	\begin{itemize}
		\item[(i)]$\eta(1,v)=v$ if $v\in \mathcal{S}(c)\setminus E^{-1}([M_c-2\vare_2,M_c+2\vare_2])$;
		\item[(ii)]$\eta(1,E^{M_c+\vare_2}\cap B_{u_c}(\delta))\subset E^{M_c-\vare_2}$;
		\item[(iii)]$E(\eta(1,v))\leq E(v)$ for all $v\in \mathcal{S}(c)$.
	\end{itemize}
	Here,  $E^\kappa:=\{v\in \mathcal{S}(c):E(v)\leq \kappa\}$. For $(i)$, by definition we see that $f(v)=0$, thus $\frac{d}{dt}\eta(t,v)|_{t=0}=W(v)=0$ and $\eta(t,v)\equiv \eta(0,v)=v$. For $(iii)$, using \eqref{eq:5.9} and the non-negativity of $g$ we obtain
	\begin{align*}
		E(\eta(1,v))&=E(v)+\int_0^1\frac{d}{ds}E(\eta(s,v))\,ds\nonumber\\
		&=E(v)-\int_{s\in[0,1],\eta(s,v)\in \tilde{S}(c)}\la E'(\eta(s,v)),f(\eta(s,v))\|Y(\eta(s,v))\|^{-1}_{H_{x,y}^1}Y(\eta(s,v))\ra\,ds\nonumber\\
		&\leq E(v)-\frac12\int_{s\in[0,1],\eta(s,v)\in \tilde{S}(c)}f(\eta(s,v))\|E'|_{\mathcal{S}(c)}(\eta(s,v))\|_*\,ds\leq E(v).
	\end{align*}
For $(ii)$, we first show that for $v\in E^{M_c+\vare_2}\cap B_{u_c}(\delta)$ one has $\eta(t,v)\in B_{u_c}(2\delta)$ for all $t\in[0,1]$. This follows from
	\begin{align*}
		\|\eta(t,v)-v\|_{H_{x,y}^1}=\|\int_{s\in[0,1],\eta(s,v)\in \tilde{S}(c)} f(\eta(t,v))\|Y(\eta(t,v))\|^{-1}_{H_{x,y}^1}Y(\eta(t,v))\,ds\|_{H_{x,y}^1}
		\leq t\delta\leq \delta.
	\end{align*}
	By \eqref{eq:5.10} this implies particularly that $\|E'|_{\mathcal{S}(c)}(\eta(t,v))\|_*\geq \sigma$. Consequently, using \eqref{eq:5.9}, $0\leq f\leq \delta$ and $\vare_2\leq \sigma\delta/4$ we obtain
	\begin{align*}
		E(\eta(1,v))
		&\leq E(v)-\frac12\int_{s\in[0,1],\eta(s,v)\in \tilde{S}(c)}f(\eta(s,v))\|E'|_{\mathcal{S}(c)}(\eta(s,v))\|_*\,ds\nonumber\\
		&\leq M_c+\vare_2-\frac{\sigma\delta}{2}\leq M_c-\vare_2.
	\end{align*}
	Next, we recall the function $g$ defined by \eqref{eq:5.8} by setting $u=u_c$ therein. We claim that there exist $0<t_0\ll 1$ and $t_1\gg 1$ such that $g\in\Gamma(c)$. Indeed, from the proof of Lemma \ref{lem:5.5} it suffices to show $2a\|\pt_y u\|_2^2+b\|\pt_y u\|_2^4 \leq 4(M_c-\vare)$. Define the scaling operator $T_\kappa$ by
	\begin{align}\label{eq:5.11}
		T_\kappa u(x,y):=\kappa^{\frac{2}{\alpha}}u(\kappa x,y).
	\end{align}
	Then
	\begin{align*}
		\|\pt_x (T_\kappa u)\|_2^2=\kappa^{1+\frac{4}{\alpha}}\|\pt_x u\|_2^2,\quad
		\|T_\kappa u\|_{\alpha+2}^{\alpha+2}=\kappa^{1+\frac{4}{\alpha}}\|u\|_{\alpha+2}^{\alpha+2},\\
	\end{align*}
	and
	\begin{align*}
		\|\pt_y (T_\kappa u)\|_2^2=\kappa^{\frac{4}{\alpha}-1}\|\pt_y u\|_2^2,\quad
		\|T_\kappa u\|_2^2=\kappa^{\frac{4}{\alpha}-1}\|u\|_2^2.
	\end{align*}	
	Combining with Lemma \ref{lem:5.3} and the fact that $u_c$ is a minimizer of $M_c$ we infer that $\frac{d}{d\kappa}E(T_\kappa u_c)|_{\kappa=1}\geq 0$, or  equivalently
	\begin{align}\label{eq:5.12}
		\frac{\alpha+4}{2(\alpha+2)}\|u_c\|_{\alpha+2}^{\alpha+2}&\geq a\|\pt_y u_c\|_{2}^2+b\|\pt_y u_c\|_{2}^4+b\|\pt_x u_c\|_{2}^2\|\pt_y u_c\|_{2}^2.
	\end{align}
	Hence
	\begin{align*}
		M_c=E(u_c)-\frac{1}{4}P(u_c)&\geq\frac{a}{2}\|\pt_y u_c\|_{2}^2+\frac{b}{4}\|\pt_y u_c\|_{2}^4+\frac{b}{4}\|\pt_x u_c\|_{2}^2\|\pt_y u_c\|_{2}^2+\frac{\alpha-8}{8(\alpha+2)}\|u_c\|^{\alpha+2}_{\alpha+2}	\\
		&\geq\frac{3\alpha}{8(\alpha+4)}\left( 2a\|\pt_y u_c\|_{2}^2+b\|\pt_y u_c\|_{2}^4\right) ,
	\end{align*}
	then we obtain
	\begin{align*}
		2a\|\pt_y u_c\|_{2}^2+b\|\pt_y u_c\|_{2}^4\leq 4M_c\left( 1-\frac{4\alpha-32}{12\alpha}\right) \leq 4M_c(1-\vare),
	\end{align*}
	as desired. By (i) and Lemma \ref{lem:5.5} we know that $\eta(1,g(t))\in \Gamma(c)$. We shall finally prove
	\begin{align*}
		\max_{t\in[0,1]}E(\eta(1,g(t)))<M_c,
	\end{align*}
	which contradicts the characterization $M_c=\gamma_c$ deduced from Lemma \ref{lem:5.5} and closes the desired proof. Notice by definition of $g$ and Lemma \ref{lem:5.2} we have $E(g(t))\leq E(u_c)=M_c$ for all $t\in[0,1]$. Thus only the following scenarios can happen:
	\begin{itemize}
		\item[(a)] $g(t)\in \mathcal{S}(c)\setminus B_{u_c}(\delta)$. By (iii) and Lemma \ref{lem:5.2} (5) we have
		\begin{align*}
			E(\eta(1,g(t)))\leq E(g(t))<E(u_c)=M_c.
		\end{align*}
		
		\item[(b)]$g(t)\in E^{m_c-\vare_2} $. By (iii) we have
		\begin{align*}
			E(\eta(1,g(t)))\leq E(g(t))\leq M_c-\vare_2<M_c.
		\end{align*}
		
		\item[(c)]$g(t)\in E^{-1}([m_c-\vare_2,m_c+\vare_2])\cap B_{u_c}(\delta) $. By (ii) we have
		\begin{align*}
			E(\eta(1,g(t)))\leq  M_c-\vare_2<M_c.
		\end{align*}
	\end{itemize}
	This completes the desired proof.
\end{proof}
Next, we give the proof of Theorem \ref{th:1.3}, and using a much handier characterization for $M_c$ due to Le Coz \cite{Le2008ANS} that is more useful for our analysis.
\begin{proof}[Proof of Theorem \ref{th:1.3}]
	We split our proof into four steps.
	
	\textit{\textbf{Step 1:} $M_c=\widetilde{M}_c$.}
	Here
	\begin{align*}
		\widetilde{M}_c:=\inf\{\mI(u):u\in \mathcal{S}(c),P(u)\leq 0\},
	\end{align*}
	and 
	$$
	\mI(u):=E(u)-\frac{1}{4}P(u)=\frac{a}{2}\|\pt_y u\|_{2}^2+\frac{b}{4}\|\pt_y u\|_{2}^4+\frac{a}{4}\|\pt_x u\|_{2}^2+\frac{b}{4}\|\pt_x u\|_{2}^2\|\pt_y u\|_{2}^2+\frac{\alpha-8}{8(\alpha+2)}\|u\|^{\alpha+2}_{\alpha+2}.
	$$ 
	We now prove $M_c=\widetilde{M}_c$. Since $V(c)\subset\{u:u\in \mathcal{S}(c),P(u)\leq 0\}$, 	\begin{align*}
		\widetilde{M}_c
		\leq\inf\{\mI(u):u\in V(c)\}
		=\inf\{E(u):u\in V(c)\}=M_c,
	\end{align*}
	On the other hand, let $\{u_k\}\subset \mathcal{S}(c)$ be a minimizing sequence for the variational problem $\widetilde{M}_c$, i.e.
	\begin{align}\label{eq:5.13}
		\mI(u_k)=\widetilde{M}_c+o_k(1),\quad
		P(u_k)\leq 0\quad\forall\,n\in\N.
	\end{align}
	By Lemma \ref{lem:5.2} we know that there exists some $t_k\in(0,1]$ such that $P(u_k\ast{t_k})=0$. Thus
	\begin{align*}
		M_c\leq E(u_k\ast{t_k})=\mI(u_k\ast{t_k})\leq \mI(u_k)=\widetilde{M}_c+o_k(1).
	\end{align*}
	Sending $k\to\infty$ we infer that $M_c\leq \widetilde{M}_c$. We complete the proof of Step 1.

	\textit{\textbf{Step 2:} Existence of a non-negative  solution of problem \eqref{eq:1.1}.}
	  Let $\{u_k\}\subset \mathcal{S}(c)$ be a minimizing sequence satisfying \eqref{eq:5.13}. Since $E(|u|)=E(u)$ for $u\in \mathcal{S}(c)$,  we may  assume that $u_k\geq 0$. By Lemma \ref{lem:5.1} and Step 1 we know that $\widetilde{M}_c<\infty$, hence
	\begin{align*}
		\infty>\widetilde{M}_c+o_k(1)=I(u_k)=\frac{a}{2}\|\pt_y u_k\|_{2}^2+\frac{b}{4}\|\pt_y u_k\|_{2}^4+\frac{a}{4}\|\pt_x u_k\|_{2}^2+\frac{b}{4}\|\pt_x u_k\|_{2}^2\|\pt_y u_k\|_{2}^2+\frac{\alpha-8}{8(\alpha+2)}\|u_k\|^{\alpha+2}_{\alpha+2}.
	\end{align*}
	Then we conclude that $\{u_k\}$ is a bounded sequence in $H_{x,y}^1$. Now using \eqref{eq:5.1} and Lemma \ref{lem:non-vanishing} we may find some non-negative $u\in H_{x,y}^1\setminus\{0\}$ such that $u_k\rightharpoonup u$ weakly in $H_{x,y}^1$. By weakly lower semicontinuity of norms we deduce
	\begin{align}
		\|u\|_2^2=:c_1\in(0,c]\quad \text{and}\quad\mI(u)\leq \widetilde{M}_c.\label{eq:5.14}
	\end{align}
	We next show $P(u)\leq 0$. Assume the contrary $P(u)>0$. By Brezis-Lieb lemma \cite{willem1997} and \cite{Chen2017CV},  we infer that
	\begin{align*}
		\|u_k-u\|_2^2&=c-c_1+o_k(1),
	\end{align*}
	and
	\begin{align*}
	P(u_k)=&	P(u_k-u)+P(u)\\
	&+2b\|\pt_x (u_k-u)\|_{2}^2\|\pt_x u\|_{2}^2+b\|\pt_x (u_k-u)\|_{2}^2\|\pt_y u\|_{2}^2+b\|\pt_x u\|_{2}^2\|\pt_y (u_k-u)\|_{2}^2+o_k(1).
	\end{align*}
	Since $P(u_k)\leq 0$, we have
	\begin{align*}
		P(u_k-u)&\leq -P(u)+o_k(1).
	\end{align*}
	Therefore, for all sufficiently large $k$ we know that $\|u_k-u\|_2^2\in[0,c)$ and $P(u_k-u)<0$. By Lemma \ref{lem:5.2} we know that there exists some $t_k\in(0,1)$ such that $P((u_k-u)\ast{t_k})=0$. Consequently, Lemma \ref{lem:5.3}, Brezis-Lieb lemma and Step 1 yield
	\begin{align*}
		\widetilde{M}_{c}&\leq \mI((u_k-u)\ast{t_k})\leq\mI(u_k-u)\\
		&=\mI(u_k)-\mI(u)-\frac{b}{2}\|\pt_y (u_k-u)\|_{2}^2\|\pt_y u\|_{2}^2-\frac{b}{2}\|\pt_x u\|_{2}^2\|\pt_y (u_k-u)\|_{2}^2-\frac{b}{2}\|\pt_x (u_k-u)\|_{2}^2\|\pt_y u\|_{2}^2+o_k(1)\\
		&=\widetilde{M}_{c}-\mI(u)+o_k(1).
	\end{align*}
	Sending $k\to\infty$ and using the non-negativity of $\mI(u)$ we obtain $\mI(u)=0$,  which is a contradiction to $u\neq0$ and thus $P(u)\leq 0$. If $P(u)<0$, then again by Lemma \ref{lem:5.2} we find some $s\in(0,1)$ such that $P(u\ast s)=0$. But then using Lemma \ref{lem:5.3}, Step 1 and the fact $c_1\leq c$
	\begin{align*}
		\widetilde{M}_{c_1}\leq \mI(u\ast s)<\mI(u)\leq \widetilde{M}_c\leq \widetilde{M}_{c_1},
	\end{align*}
	a contradiction. We conclude therefore $P(u)=0$.	
	Thus $u$ is a minimizer of $M_{c_1}$. From Lemma \ref{lem:5.6} we know that $u$ is a solution of problem \eqref{eq:1.1}.
	
	\textit{\textbf{Step 3:} Positivity of $\omega$.}
	First we prove that $\omega$ is non-negative. Testing problem \eqref{eq:1.1} with $u$ and followed by eliminating $a\|\pt_x u\|_2^2+b\|\pt_x u\|_2^4$ using $P(u)=0$ we obtain
	\begin{align}
		a\|\pt_y u\|_2^2+b\|\pt_y u\|_2^4+b\|\pt_x u\|_2^2\|\pt_y u\|_2^2+\omega\|u\|_2^2=\frac{\alpha +4}{2(\alpha+2)}\|u\|_{\alpha+2}^{\alpha+2}.\label{eq:5.15}
	\end{align}
	Combining with \eqref{eq:5.12} we infer that $\omega\|u\|_2^2\geq 0$. Since $u\neq 0$ we conclude $\omega\geq 0$. We next assume the contrary $\omega=0$. Therefore that $u$ satisfies the equation
	\begin{align}\label{eq:5.16}
		-\left( a+bB\right) \Delta_{x,y}u=|u|^\alpha u,
	\end{align}
	where $B=\lim_{k\to\infty}\|\nabla_{x,y} u_k\|_{2}^2$.
	For $n\in\N$, let $\phi_k\in C_0^\infty(\R;[0,1])$ be a radially symmetric decreasing function such that $\phi_k(t)\equiv 1$ on $|t|\leq k$, $\phi_k(t)\equiv 0$ for $|t|\geq k+1$ and $\sup_{k\in\N}\|\phi_k\|_{L^\infty}\leq C$. Define 
	$$ D_k:=\{x\in\R:|x|\in(k,k+1)\}.$$
	Since $u\neq 0$ is non-negative, by monotone convergence theorem we know that
	\begin{align*}
		\int_{\R\times\T}u^{\alpha+1}\phi_k\,dxdy>0\quad\text{for all}\quad k\gg 1.
	\end{align*}
	On the other hand, using H\"older and the fact that $\mathrm{supp}\,\pt_x\phi_k\subset D_k$,  we see that
	\begin{align}
		\int_{\R\times \T}\pt_x u\cdot\pt_x \phi_k\,dxdy
		&\leq \|\pt_x u\|_{L^{2}(D_k\times\T)}\|\pt_x \phi_k\|_{L^2(D_k\times\T)}\nonumber\\
		&\leq  C\|\pt_x u\|_{L^{2}(D_k\times\T)}.\label{contra1}
	\end{align}
	Then
	$\infty>\|\pt_x u\|_{2}^{2}\geq \sum_{k\geq 1}\|\pt_x u\|^{2}_{L^{2}(D_k\times\T)}$
	yields $\|\pt_x u\|_{L^{2}(D_k\times\T)}=o_k(1)$. By the fact that $\phi_k$ is independent of $y$ we also know $\int_{\R\times\T}(-\pt_y^2 u)\phi_k\,dxdy=0$.
	Summing up, by testing \eqref{eq:5.16} with $\phi_k$ and rearranging terms we obtain
	\begin{align*}
		o_k(1)=\int_{\R\times\T}u^{\alpha+1}\phi_k\,dxdy>0
	\end{align*}
	for $k$ sufficiently large, a contradiction. This completes the proof of Step 3.
	
\textit{	\textbf{Step 4:} $c_1=c$.}
	Assume that $c_1<c$. By Lemma \ref{lem:5.3} and \eqref{eq:5.14} we know that $M_{c_1}$ is a local minimizer of the mapping $c\mapsto M_c$, which in turn implies that the inequality in \eqref{eq:5.12} is in fact an equality. Now, combine with  \eqref{eq:5.15} we infer that $\omega\|u\|_2^2=0$, which is a contradiction since $\omega>0$ and $u\neq 0$. We thus conclude $\|u\|_2^2=c$. 
	That $u$ is positive follows immediately from the strong maximum principle. This completes the desired proof.
\end{proof}

\subsection{$y$-dependence of normalized ground states}\label{sub2.2}
In this subsection we prove Theorem \ref{th:1.4}. Similar as the proof of Theorem \ref{th:1.2},
we consider the variational problems $M_{1,\lambda}$ and $\nu_{1,\tau}$ defined in \eqref{M c lambda} and \eqref{nu c tau} respectively.

\begin{lemma}\label{lem:6.1}
	Let  $\alpha\in(8,\infty)$ and $u_\tau$ be a minimizer of $\nu_{1,\tau}$ which also satisfies
	\begin{align}
		&-\left( a+\frac{b}{2\tau^{\frac{4}{\alpha}-1}}\|\nabla_{y} u_{\tau}\|^2_{2}+\frac{b}{\tau^{\frac{4}{\alpha}+1}} \|\pt_xu_\tau\|_2^2\right) \Delta_{x}u_\tau\nonumber\\
		&-\left( a\tau^2+\frac{b}{2\tau^{\frac{4}{\alpha}-1}}\|\pt_x u_{\tau}\|^2_{2}+b\tau^{3-\frac{4}{\alpha}}\|\nabla_{y}u_\tau\|^2_2\right) \pt^2_{y}u_\tau +\omega_\tau w_\tau=|u_\tau|^\alpha u_\tau\quad\text{on $\R\times \T$}\label{eq:6.1}
	\end{align}
	for some $\omega_\tau>0$ (whose existence is guaranteed by Theorem \ref{th:1.3}). Then we have
	\begin{align}
		\lim_{\tau\to\infty}\nu_{1,\tau}=(2\pi)^n\widehat{\nu}_{(2\pi)^{-n},\infty}\label{eq:6.2}
	\end{align}
	(here $\widehat{\nu}_{c,\infty}$ is defined in \eqref{hat nu infty}) and
	\begin{align}
		&\lim_{\tau \to \infty}\tau^2\|\pt_{y} u_{\tau}\|^2_{2}=0,\label{eq6.3}\\ &\lim_{\tau \to \infty}\tau^{3-\frac{4}{\alpha}}\|\pt_{y} u_{\tau}\|^4_{2}=0,\\
		&\lim_{\tau \to \infty}\frac{1}{\tau^{\frac{4}{\alpha}-1}}\|\pt_x u_{\tau}\|^2_{2}\|\pt_{y} u_{\tau}\|^2_{2}=0,\\ &\lim_{\tau \to \infty}\frac{1}{\tau^{\frac{4}{\alpha}+1}}\|\pt_x u_{\tau}\|^4_{2}=0.\label{eq:6.3}
	\end{align}
\end{lemma}
\begin{proof}
	By assuming that  $u\in V(1)$ is independent of $y$ and $\widehat{E}^{\infty}(u) = \widehat{\nu}_{(2\pi)^{-1},\infty}$.  Thus we get $u\in V(1)$ and
	\begin{align*}
		\nu_{1,\tau} \leq E^{\tau}(u) &= 2\pi\left(\frac{a}{2}\|\pt_x u\|^2_{L_x^2}+\frac{2\pi b}{4\tau^{\frac{4}{\alpha}+1}}\|\pt_x u\|^4_{L_x^2} -\frac{1}{\alpha+2}\|u\|^{\alpha+2}_{L_x^{\alpha+2}}\right) \\
		&= 2\pi\left(\widehat{E}^{\infty}(u)+\frac{2\pi b}{4\tau^{\frac{4}{\alpha}+1}}\|\pt_x u\|^4_{L_x^2} \right)=(2\pi)^n\widehat{\nu}_{(2\pi)^{-1},\infty} + o_\tau(1),
	\end{align*}
	where $\lim_{\tau \to \infty}o_\tau(1) = 0$, then
	\begin{equation}\label{eq:6.4}
		\limsup_{\tau \to \infty}\nu_{1,\tau} \leq (2\pi)^n\widehat{\nu}_{(2\pi)^{-1},\infty}.
	\end{equation}
	Since $u_\tau$ be a minimizer of $\nu_{1,\tau}$ and $\alpha>8$, then we have 
	\begin{align*}
		P^\tau(u_\tau):=a\|\pt_x u_\tau\|_{2}^2+\dfrac{b}{\tau^{\frac{4}{\alpha}-1}}\|\nabla_{x,y} u_\tau\|_{2}^2\|\pt_x u_\tau\|_{2}^2+\dfrac{b}{\tau^{\frac{4}{\alpha}+1}}\|\pt_x u_\tau\|_{2}^4-\frac{\alpha}{2(\alpha+2)}\| u_\tau\|_{\alpha+2}^{\alpha+2}=0,
	\end{align*}
	(here $P^\tau(u)$ is defined in \eqref{pohozave tau}) and
	\begin{align}
		\infty>	\nu_{1,\tau}=&E^\tau(u_\tau)-\frac{2}{\alpha}P^\tau(u_\tau)\nonumber\\
		=&\frac{a\tau^2}{2}\|\pt_y u_\tau\|_2^2+\frac{b\tau^{3-\frac{4}{\alpha}}}{4}\|\pt_y u_\tau\|_2^4+\frac{b}{\tau^{\frac{4}{\alpha}-1}}\left( \frac{1}{2}-\frac{2}{\alpha}\right)\|\pt_x u_{\tau}\|^2_{2}\|\nabla_{y} u_{\tau}\|^2_{2}\nonumber\\
		&+a\left( \frac{1}{2}-\frac{2}{\alpha}\right) \|\pt_x u_\tau\|_2^2+\frac{b}{\tau^{\frac{4}{\alpha}+1}}\left( \frac{1}{4}-\frac{2}{\alpha}\right) \|\pt_x u_\tau\|_2^4\nonumber\\
		\geq& \frac{a\tau^2}{2}\|\pt_y u_\tau\|_2^2+\frac{b\tau^{3-\frac{4}{\alpha}}}{4}\|\pt_y u_\tau\|_2^4\label{eq:6.5}
	\end{align}
	showing that
	\begin{align}\label{eq:6.6}
		\lim_{\tau\to\infty}\|\pt_{y} u_{\tau}\|^2_{2}=0.
	\end{align}
	Using \eqref{eq:6.4} and \eqref{eq:6.5} we infer that
	\begin{align*}
		\|\pt_x u_\tau\|_2^2\leq C\nu_{1,\tau}\leq 2\pi C\widehat{\nu}_{(2\pi)^{-1},\infty}<\infty.
	\end{align*}
	Therefore $\{u_\tau\}$ is a bounded sequence in $H_{x,y}^1$, whose weak limit is denoted by $u$. Using \eqref{eq:5.1} and Lemma \ref{lem:non-vanishing} we may also assume that $u\neq 0$. By \eqref{eq:6.6} we know that $u$ is independent of $y$ and thus $u\in H_x^1$. Moreover, using weakly lower semicontinuity of norms we know that $\|u\|_{L^2_x}^2\in(0,(2\pi)^{-1}]$. On the other hand, using $P^\tau(u_\tau)=0$, \eqref{eq:6.1} and $u_\tau\in \mathcal{S}(1)$ we obtain
	\begin{align*}
		0<\omega_\tau=\frac{2\alpha +(4-\alpha)}{2(\alpha+2)}\|u_\tau\|_{\alpha+2}^{\alpha+2}-a\tau^2\|\pt_y u_\tau\|_2^2-b\tau^{3-\frac{4}{\alpha}}\|\pt_y u_\tau\|_2^4-\dfrac{b}{\tau^{\frac{4}{\alpha}-1}}\|\nabla_{x,y} u_\tau\|_{2}^2\|\pt_x u_\tau\|_{2}^2\leq C\|u_\tau\|_{\alpha+2}^{\alpha+2}.
	\end{align*}
	Then $\{\omega_\tau\}$ is a bounded sequence in $(0,\infty)$ and  $\omega_\tau\to\omega$ as $\tau\to\infty$. We now test \eqref{eq:6.1} with $\phi\in C_0^\infty(\mathbb{R})$ and integrate both sides over $\mathbb{R}\times\T$. Notice particularly that the term $\int_{\mathbb{R}\times\T}  \pt_y^2 u_\tau \phi\,dxdy=0$ for any $\tau>0$ since $\phi$ is independent of $y$. Using  $u_\tau\rightharpoonup u$ in $H_{x,y}^1$, by sending $\tau\to\infty$ we obtain
	\begin{align*}
		-a\Delta_x u+\omega u=|u|^\alpha u\quad\text{in $\mathbb{R}$}.
	\end{align*}
	In particular, by Porposition \ref{lem:1.5} we know that $\widehat{P}^\infty(u)=0$ and $\omega>0$. Combining with weakly lower semicontinuity  we deduce
	$$
	2\pi\widehat{E}^\infty(u)=2\pi\widehat{I}^\infty(u)\leq\liminf_{\tau\to\infty}\mI_\tau(u_\tau)= \liminf_{\tau\to\infty}E_\tau(u_\tau)\leq 2\pi \widehat{\nu}_{(2\pi)^{-1},\infty},
	$$
	where $I^\tau(u):=E^\tau(u)-\frac{1}{4}P^\tau(u)$ and $\widehat{I}^\infty(u)=\widehat{E}^\infty(u)-\frac{1}{4}\widehat{P}^\infty(u)$ (here $\widehat{P}^\infty(u)$ is defined in \eqref{hat pohozave infty}). However, by Proposition \ref{lem:1.5} the mapping $c\mapsto \widehat{\nu}_{c,\infty}$ is strictly monotone decreasing on $(0,\infty)$, from which we conclude that $u\in\mathcal{S}((2\pi)^{-1})$ and $u$ is a minimizer of $\widehat{\nu}_{(2\pi)^{-1},\infty}$. 
	Finally, using the weakly lower semicontinuity of norms we conclude
	\begin{align}
		\nu_{1,\tau}=E^\tau(u_\tau)-\frac{2}{\alpha}P(u_\tau)=&\frac{a\tau^2}{2}\|\pt_y u_\tau\|_2^2+\frac{b}{2\tau^{\frac{4}{\alpha}-1}}\|\pt_x u_{\tau}\|^2_{2}\|\nabla_{y} u_{\tau}\|^2_{2}+\frac{b\tau^{3-\frac{4}{\alpha}}}{4}\|\pt_y u_\tau\|_2^4\nonumber\\
		&+a\left( \frac{1}{2}-\frac{2}{\alpha}\right) \|\pt_x u_k\|_2^2+\frac{b}{\tau^{\frac{4}{\alpha}+1}}\left( \frac{1}{4}-\frac{2}{\alpha}\right) \|\pt_x u_k\|_2^4\nonumber\\
		\geq& a\left( \frac{1}{2}-\frac{2}{\alpha}\right) \|\pt_x u_\tau\|_2^2\nonumber\\
		\geq& a\left( \frac{1}{2}-\frac{2}{\alpha}\right) 2\pi\|\pt_x u\|_{L_x^2}^2+o_\tau(1)\nonumber\\
		=&2\pi \wmH^\infty(u)+o_\tau(1)\geq 2\pi\widehat{\nu}_{(2\pi)^{-1},\infty}+o_\tau(1).\label{eq:6.7}
	\end{align}
	Then \eqref{eq:6.2} follows now from \eqref{eq:6.4} and \eqref{eq:6.7}. Finally, by keeping 
	$$
	\frac{a\tau^2}{2}\|\pt_{y} u_\tau\|^2_{2}+\frac{b}{2\tau^{\frac{4}{\alpha}-1}}\|\pt_x u_{\tau}\|^2_{2}\|\pt_{y} u_{\tau}\|^2_{2}+\frac{b\tau^{3-\frac{4}{\alpha}}}{4}\|\pt_{y} u_\tau\|^4_{2}
	+\frac{b}{\tau^{\frac{4}{\alpha}+1}}\left( \frac{1}{4}-\frac{2}{\alpha}\right)\|\pt_x u_\tau\|^4_{2}
	$$ 
	in \eqref{eq:6.7} we have
	\begin{align}
		\nu_{1,\tau}\geq &\frac{a\tau^2}{2}\|\pt_{y} u_\tau\|^2_{2}+\frac{b}{2\tau^{\frac{4}{\alpha}-1}}\|\pt_x u_{\tau}\|^2_{2}\|\pt_{y} u_{\tau}\|^2_{2}+\frac{b\tau^{3-\frac{4}{\alpha}}}{4}\|\pt_{y} u_\tau\|^4_{2}
		+\frac{b}{\tau^{\frac{4}{\alpha}+1}}\left( \frac{1}{4}-\frac{2}{\alpha}\right)\|\pt_x u_\tau\|^4_{2}\nonumber\\
		&+2\pi\widehat{\nu}_{(2\pi)^{-1},\infty}+o_\tau(1).\label{eq:6.8}
	\end{align}
	Thus \eqref{eq6.3}-\eqref{eq:6.3} follow now from \eqref{eq:6.2} and taking $\tau$ in \eqref{eq:6.8} to infinity. This completes the desired proof.
\end{proof}

\begin{lemma}\label{lem:6.2}
	Let $u_\tau, u$ be the functions given in Lemma \ref{lem:6.1}, then $u_\tau\to u$ in $H_{x,y}^1$ as $\tau\to\infty$, and $\widehat{E}^{\infty}(u)=\widehat{\nu}_{(2\pi)^{-1},\infty}$.
\end{lemma}
\begin{proof}
	In the proof of Lemma \ref{lem:6.1} we see that all the inequalities involving the weakly lower semicontinuity of norms are in fact equalities. This particularly implies $\|u_\tau\|_{H_{x,y}^1}\to \|u\|_{H_{x,y}^1}=2\pi\|u\|_{H_x^1}$ as $\tau\to\infty$, which in turn implies the strong convergence of $u_\tau\to u$ in $H_{x,y}^1$.
\end{proof}

\begin{lemma}\label{lem:6.3}
	Let  $\alpha\in(8,\infty)$. Then there exists some $\tau^*>0$ such that for all $\tau>\tau^*$ we have $\pt_y u_\tau=0$.
\end{lemma}

\begin{proof}
	This can be proved similarly by using the arguments applied in the proofs of Lemmas \ref{lem:4.5}, \ref{lem:6.1} and \ref{lem:6.2}, we thus omit the details here.
\end{proof}

\begin{proof}[Proof of Theorem \ref{th:1.4}]
	\textbf{Case (i)}:
	Recall the quantitiy $\widehat{M}_{c,\lambda}$ defined by \eqref{M c lambda}.  
	Define $\rho$ as \eqref{eq rho}.	We chose the value of $\beta$ by
	$$
	\beta=(\pi-\delta)^{-1}\left( \frac{\alpha+3}{3}\right) ^{\frac{1}{\alpha}},
	$$ 
	such that 
	$$
	\|\rho\|_{L^{\alpha+2}_y}^{\alpha+2}=\|\rho\|_{L_y^2}^2= \frac{2}{3}\left( \frac{\alpha+3}{3}\right) ^{\frac{2}{\alpha}}(\pi-\delta).
	$$
	Let $\delta\in(0,\pi)$ such that $\delta=\pi-\frac{\pi}{2}\bg(\frac{3}{\alpha+3}\bg)^{\frac2\alpha}$. 	By direct calculation we see that $\rho\in H_y^1$ and
	\begin{align*}
		2\pi>1=\|\rho\|_{L_y^2}^2=\|\rho\|_{L_y^{\alpha+2}}^{\alpha+2}.
	\end{align*}
		Using Theorem \ref{thm:2.4} we know that $\widehat{M}_{\|\rho\|^{-2n}_{L^2(\mathbb{T})},\lambda}>0$  has a minimizer $Q_\ld$, then
		\begin{align*}
			\widehat{P}_\lambda(Q_\lambda)=a\lambda^{\frac{3}{\alpha-2}+\frac{1}{2}}\|\pt_x Q_\lambda\|_{L_x^2}^2+b\|\pt_x Q_\lambda\|_{L_x^2}^4-\frac{\alpha}{2(\alpha+2)}\|Q_\lambda\|_{L_x^{\alpha+2}}^{\alpha+2}=0.
		\end{align*}
		  Arguing as in Lemma \ref{lem:5.1} we then infer that $\{Q_\lambda\}$ for $\lambda\in[0,1)$ is a bounded sequence in $H^1_{x}$.
	Next, we show that $M_{1,\ld}<(2\pi)^{n} \widehat{M}_{(2\pi)^{-n},\lambda}$ for all $\ld\ll 1$.  
	Next,we define the function $\psi$ by
	\begin{align*}
		\psi_\lambda(x,y)&:=Q_\ld(x) \rho(y).
	\end{align*}
	 We conclude that $\psi\in \mathcal{S}(1)$, and for given $\lambda$ let $t_\lambda\in(0,\infty)$ be given such that $P_\lambda(\psi*t_\lambda)=0$, where $\psi*t_\lambda$ is defined by \eqref{eq:L2scaling}. Since $\|\rho\|_{L_y^2}^2=\|\rho\|_{L_y^{\alpha+2}}^{\alpha+2}$, we have $P_0(\psi_0)=0$, This, in conjunction with standard continuity arguments, implies $t_\lambda\to1$ as $\lambda\to0$.
	As in Lemma \ref{lem:5.3}, the function $c\mapsto \widehat{M}_{c,\lambda}$ is strictly decreasing on $(0,\infty)$. Then $\|\rho\|_{L_y^2}^{2}\widehat{M}_{\|\rho\|_{L_y^2}^{-2},\lambda}=\widehat{M}_{1,\lambda}
	<2\pi\widehat{M}_{(2\pi)^{-1},\lambda}$.	Consequently,
	\begin{align}
		M_{1,\ld}\leq E_\ld(\psi\ast{t_\lambda})
		=&\|\rho\|_{L^2_y}^{2}\left( \frac{a\lambda^{\frac{3}{\alpha-2}+\frac{1}{2}}}{2}\|\pt_x  (Q_\ld\ast{t_\lambda})\|_{L_x^2}^2+\frac{b}{4}\|\pt_x (Q_\ld\ast{t_\lambda})\|_{L_x^2}^4
		-\frac{1}{\alpha+2}\| (Q_\ld\ast{t_\lambda})\|_{L_x^{\alpha+2}}^{\alpha+2}\right) \nonumber\\
		&+\frac{b}{4}\|\rho\|_{L^2_y}^{2}\left( \|\rho\|_{L^2_y}^{2}-1\right) \|\pt_x (Q_\ld\ast{t_\lambda})\|_{L_x^2}^4-\frac{\|\rho\|_{L^2_y}^{2}}{\alpha+2}
		\left( \|\rho\|_{L^2_y}^{-2}\|\rho\|_{L^2_y}^{\alpha+2}-1\right) 
		\|(Q_\ld\ast{t_\lambda})\|_{L_x^{\alpha+2}}^{\alpha+2}\nonumber\\
		&+ \frac{b\lambda^2}{4\|\rho\|_{L^2_y}^{4}}\|\pt_y\rho\|_{L_y^2}^4+\frac{b\lambda}{2}\|\pt_x(Q_\ld\ast{t_\lambda})\|^2_{L^2_x}\|\pt_{y} \rho\|^2_{L^2_y}
		+\frac{a\lambda^{\frac{3}{\alpha-2}+\frac{3}{2}}}{2\|\rho\|_{L^2_y}^{2}}\|\pt_y \rho\|_{L_y^2}^2 \nonumber\\
		=&\|\rho\|_{L^2_y}^{2}\widehat{M}_{\|\rho\|_{L^2_y}^{-2},\lambda}+o_\lambda(1)	<2\pi\widehat{M}_{(2\pi)^{-1},\lambda}+o_\lambda(1).\nonumber
	\end{align}
Then we infer that 
$$
\lim_{\ld\to 0}M_{1,\lambda}<2\pi\widehat{M}_{(2\pi)^{-1},\lambda},
$$
	and there exists $\lambda^*>0$ such that $M_{1,\ld}<2\pi\widehat{M}_{(2\pi)^{-1},\lambda}$ for $\ld<\ld^*$. That a minimizer $u_\ld$ of $M_{1,\ld}$ has non-trivial $y$-dependence for $\ld<\ld^*$, this is $ \pt_{y} u_\ld\neq0$.

Similar as the proof of Theorem \ref{th:1.2} $(i)$, for $c>0$ and let $\theta_c:=c^{\frac{\alpha-2}{\alpha-8}}$, direct calculation  shows that for $\lambda=\theta_c^{2}$ we obtain $u\mapsto T_1^{\theta_c} u$ defines a bijection between $V(c)$ and $V_\lambda(1)$, and
$$
E(u)=c^{\frac{2\alpha+8}{\alpha-8}}E_\lambda(T_1^{\theta_c} u).
$$
Thus $
M_{c}=c^{\frac{2\alpha+8}{\alpha-8}}M_{1,\lambda}.
$ with $\lambda=\theta_c^2$.
By same rescaling arguments we also infer that 
\begin{equation}
	\widehat{M}_{(2\pi)^{-n}c}=c^{\frac{2\alpha+8}{\alpha-8}}\widehat{M}_{(2\pi)^{-n},\lambda} ~~\text{with $\lambda=\theta_c^2$}.
\end{equation} 
 Notice also that the mapping $c\mapsto\theta_c$ is strictly monotone increasing on $(0,\infty)$. Thus we know that there exists some $c^{*}\in(0,\infty)$  such that for any $c\in(0,c^{*})$ we have
  $
  M_{c}<2\pi\widehat{M}_{(2\pi)^{-n}c},
  $ 
  and then for $u_c$ be a minimizer of $M_c$ we infer that  $\partial_y u_c\neq 0$ with $c\in(0,c^{*})$.

	\textbf{Case (ii)}:
	Lemma \ref{lem:6.3} yields the existence of $\tau^* \in (0,\infty)$ such that for all $\tau \in (\tau^*,\infty)$ any minimizer $u_\tau$ of $M_{1,\tau}$ satisfies $\pt_y u_\tau = 0$. 	
	This	implies that for all $\tau \in (\tau^\ast,\infty )$ we have $\nu_{1,\tau} = 2\pi \widehat{\nu}_{(2\pi )^{-1},\tau}$. 
	Similar as the proof of Theorem \ref{th:1.2} $(ii)$, direct calculation also shows that for $\tau=\kappa_c:=c^{\frac{\alpha}{\alpha-4}}$ we obtain $u\mapsto T_2^{\kappa_c} u$ defines a bijection between $V(c)$ and $V^\tau(1)$, and
	$$
	E(u)=c^{\frac{4/\alpha+1}{4/\alpha-1}}E^{\tau}(T_2^{\kappa_c} u)\quad\text{for}\quad\tau = \kappa_c,
	$$
	thus $M_{c}=c^{\frac{4/\alpha+1}{4/\alpha-1}}\nu_{1,\tau}$ for $\tau = \kappa_c$. 
	And we also infer that $\widehat{M}_{2\pi c}=\widehat{\nu}_{2\pi c,1}=c^{\frac{4/\alpha+1}{4/\alpha-1}}\widehat{\nu}_{2\pi,\tau}$
	for $\tau = \kappa_c$.  
	Notice also that the mapping $c\mapsto \kappa_c$ is strictly monotone increasing on $(0,\infty)$. Thus, if $u_c$ be a minimizer of $M_c$, we know that there exists some $c^{**}\in(0,\infty)$  such that  
for any $c\in(c^{**},\infty)$ it holds $\pt_y u_c= 0$. 
\end{proof}

{\small \noindent \text{Conflict of interest:} The authors declare that they have no conflict of interest.}

{\small \noindent \text{Data availability:} Data sharing is not applicable to this article as no datasets were generated or analysed during the current study.}


\providecommand{\href}[2]{#2}
\providecommand{\arxiv}[1]{\href{http://arxiv.org/abs/#1}{arXiv:#1}}
\providecommand{\url}[1]{\texttt{#1}}
\providecommand{\urlprefix}{URL }

\end{document}